\documentstyle[11pt]{article}
\begin{document}
\newcommand{\qed}{\hphantom{.}\hfill $\Box$\medbreak}
\newcommand{\proof}{\noindent{\bf Proof \ }}
\newtheorem{Theorem}{Theorem}[section]
\newtheorem{Lemma}[Theorem]{Lemma}
\newtheorem{Corollary}[Theorem]{Corollary}
\newtheorem{Remark}[Theorem]{Remark}
\newtheorem{Example}[Theorem]{Example}
\newtheorem{Definition}[Theorem]{Definition}
\newtheorem{Construction}[Theorem]{Construction}
\newtheorem{Question}[Theorem]{Question}
\newcommand{\B}{{\cal B}}
\newcommand{\D}{{\cal D}}
\newcommand{\E}{{\cal E}}
\newcommand{\C}{{\cal C}}
\newcommand{\F}{{\cal F}}
\newcommand{\A}{{\cal A}}
\newcommand{\Hh}{{\cal H}}
\newcommand{\Pp}{{\cal P}}
\newcommand{\G}{{\cal G}}

% 2009. 3. 1
%%%%%%%%%%%%%%%%%%%%%%%%%%%%%%%%%%%%%%%%%%%%%%%%%%%%%%%%%%%%%%%%%%%%%%%%

\begin{center}
{\Large\bf A pair of disjoint $3$-GDDs of type $g^t u^1$
*\footnote{Research of Y. Chang and J. Zhou is Supported  by National
Natural Science Foundation of China under grants 10771013 and 10831002. Research
of Y. M. Chee is supported in part by the National Research
Foundation of Singapore under Research Grant NRF-CRP2-2007-03 and by
the Nanyang Technological University under Research Grant
M58110040.}}

\vskip24pt

Yanxun Chang \\ Institute of Mathematics\\ Beijing Jiaotong
University\\ Beijing 100044, P. R. China\\ {\tt yxchang@bjtu.edu.cn}

\vskip12pt
Yeow Meng Chee\\ Division of Mathematical Sciences\\ Nanyang Technological University\\
Singapore 637371, Singapore\\ {\tt ymchee@ntu.edu.sg}

\vskip12pt Junling Zhou \\ Institute of Mathematics\\ Beijing
Jiaotong University\\ Beijing 100044, P. R. China\\ {\tt
jlzhou@bjtu.edu.cn}

\end{center}

\begin{abstract}
Pairwise disjoint 3-GDDs can be used to construct some optimal
constant-weight codes. We study the existence of a pair of disjoint
3-GDDs of type $g^tu^1$ and establish that its necessary
 conditions are also sufficient.
\medskip


\end{abstract}
\section{Introduction}

Let $X$ be a finite set of $v$ elements and $K$ a set of positive
integers. A {\em group divisible design} $K$-GDD is a triple $(X,
{\cal G}, {\cal A})$ satisfying the following properties: $(1)$
$\cal G$ is a partition of $X$ into subsets (called {\em groups});
$(2)$ $\cal A$ is a set of subsets of $X$ (called {\em blocks}),
each of cardinality from $K$, such that a group and a block contain
at most one common point; $(3)$ every pair of points from distinct
groups occurs in exactly one block. If $\cal G$ contains $u_i$
groups of size $g_i$ for $1\leq i\leq s$, then we call
$g_1^{u_1}g_2^{u_2}\cdots g_s^{u_s}$ the {\em group type} (or {\em
type}) of the GDD. If $K=\{k\}$, we write $\{k\}$-GDD as $k$-GDD. A
$k$-GDD of type $t^k$ is denoted by TD$(k,t)$ and is called a {\it
transversal design.} A $K$-GDD of type $1^v$ is commonly called a
{\em pairwise balanced design}, denoted by $(v,K,1)$-PBD. When
$K=\{k\}$, a pairwise balanced design is just a Steiner system
S$(2,k,v)$.
%When
%$K=\{k\}$, a pairwise balanced design is just a Steiner system
%$S(2,k,v)$, called a {\em balanced incomplete block design} and
%denoted by $(v,k,1)$-BIBD. A $(v,3,1)$-BIBD is called a {\it Steiner
%triple system} of order $v$, briefly by  STS$(v)$.
It is well-known that an S$(2,3,v)$ exists if and only if $v\equiv
1,3$ (mod 6).

Colbourn et al. completely settle the necessary and sufficient
conditions for the existence of $3$-GDDs of type $g^t u^1$.
\begin{Lemma}
\label{2.3}$({\rm \cite{chr}})$ Let $g$, $t$, and $u$ be nonnegative
integers. There exists a $3$-GDD of type $g^tu^1$ if and only if the
following conditions are all satisfied:
\begin{enumerate}
\item[{\rm (1)}] if $g>0$, then $t\geq3$, or $t=2$ and $u=g$,
or $t=1$ and $u=0$, or $t=0;$

\item[{\rm (2)}] $u\leq g(t-1)$ or $gt=0;$

\item[{\rm (3)}] $g(t-1)+u \equiv 0$ $({\rm mod}$ $2)$ or $gt=0;$

\item[{\rm (4)}] $gt \equiv 0$ $({\rm mod}$ $2)$ or $u=0;$

\item[{\rm (5)}] $\frac{1}{2}g^2t(t-1)+gtu \equiv 0$ $({\rm mod}$ $3).$
\end{enumerate}

\end{Lemma}

%A $K$-GDD of type $1^{v-h}$$h^1$ is commonly called an {\em
%incomplete pairwise balanced design}, denoted by $(v,h;K,1)$-IPBD.
%When $K=\{k\}$, an incomplete pairwise balanced design is called an
%{\em incomplete balanced incomplete block design}, denoted by
%$(v,h;k, 1)$-IBIBD.
Let $2\notin K$. A {\em partial group divisible design} $K$-GDD is a
triple $(X, {\cal G}, {\cal A})$ satisfying conditions (1) and (2)
of the definition of a $K$-GDD and (3') every pair of points from
distinct groups occurs in at most one block.  The {\it leave} of a
partial $K$-GDD is a graph whose edges are all the pairs which
belong to distinct groups and  do not appear in any block. A $K$-GDD
can be regarded as a partial $K$-GDD with an empty leave. Suppose
that $(X,{\cal G},{\cal B})$ and  $(X,{\cal G},{\cal B}')$ are two
partial  $K$-GDDs. If ${\cal B}$ and ${\cal B}'$ have no block in
common, $(X,{\cal G},{\cal B})$ and $(X,{\cal G},{\cal B}')$ are
said to be {\it disjoint}.

 The purpose of this paper is to
determine the existence spectrum of a pair of disjoint $3$-GDDs of
type $g^t u^1$. The problem is itself
 interesting in the theory of combinatorial designs. Also we have
a motivation lying in  a close relation between disjoint 3-GDDs and
constant-weight codes. In Chee et al. \cite{chee1}, pairwise
disjoint combinatorial designs of various types, including Steiner
systems and group divisible designs, are utilized to construct
optimal $q$-ary constant-weight codes with $q>2$. In particular,  a
pair of disjoint 3-GDDs of type $1^{6t}5^1$ is proved to exist for
any positive integer $t$, which is used in constructing optimal
3-ary constant-weight codes of Hamming distance 4 and weight 3. In
\cite{chee2}, the concept of group divisible design is generalized
to a new code named group divisible code, which is shown useful in
recursive constructions for constant-weight and constant-composition
codes. One can also find applications of disjoint group divisible
designs in the determination of more optimal  constant-weight codes
(see, for example, \cite{zhangge,zhangge2}).

In order to study the existence of two disjoint 3-GDDs, we introduce
some related notions and basic facts in this section. Let $(X,{\cal
G},{\cal A})$ be a $K$-GDD. A subset of the block set ${\cal A}$ is
called a {\it parallel class} if it contains every element of $X$
exactly once. If ${\cal A}$ can be partitioned into some parallel
classes, the GDD is called {\it resolvable.} A resolvable S$(2,3,v)$
is the well-known Kirkman triple system of order $v$, denoted by
KTS$(v)$. A KTS$(v)$ exists if and only if $v\equiv 3$ (mod 6) (see
\cite{raywilson}).

%precisely $r$
%blocks in common, then $r$ is called their {\it intersection
%number}. In particular, if $r=0$, the two $K$-GDDs  are said to be
A {\it Latin square} of order $t$ (briefly by LS$(t)$) is a $t\times
t$ array in which each cell contains a single element from a
$t$-set, such that each element occurs exactly once in each row and
exactly once in each column.  Suppose that $L=(a_{ij})$ is an
LS$(t)$ defined on and indexed by a  set $T$.  If for each $i\in T$,
$a_{ii}=i$, then the Latin square is called {\it idempotent}. If for
any $i,j\in T$, $a_{ij}=a_{ji}$, then it is called {\it symmetric}.
Suppose that $L=(a_{ij})$ and $L'=(b_{ij})$ are LS$(t)$s on a set
$T$. $L$ and $L'$ are {\it orthogonal} if every element of $T\times
T$ occurs exactly once among the $t^2$ pairs $(a_{ij},b_{ij})$,
$1\leq i,j\leq t$. %A Latin square is called {\it self-orthogonal} if
%it is orthogonal to its transpose.

 A TD$(3,t)$ is often defined on $V\times I$ with groups  $V\times\{i\},i\in I$, where $|V|=t$,
 and $|I|=3$. If the TD$(3,t)$ has a
parallel class $\{\{x\}\times I:x\in V\}$, then it is called {\it
idempotent} and denoted by ITD$(3,t)$. An ITD$(3,t)$ is equivalent
to an idempotent LS$(t)$. So when $t\geq 4$, an ITD$(3,t)$ exists.
If the block set of an ITD$(3,t)$  can be partitioned into $t$
parallel classes, one of which is the idempotent one, we call it
{\it resolvable} and denote by RITD$(3,t)$. An RITD$(3,t)$, which is
equivalent to a pair of
 orthogonal LS$(t)$s, exists if and only if $t\neq 2,6$.

 Let $(X,{\cal G},{\cal B})$ and
$(X,{\cal G},{\cal B}')$ be two ITD$(3,t)$s. They are called {\it
disjoint} if ${\cal B}$ and ${\cal B}'$ have no block in common
except the common idempotent parallel class. Similarly we have the
definition of disjoint RITDs. Note that although a resolvable
TD$(3,t)$ can always be made idempotent, two disjoint RTD$(3,t)$s do
not always mean two disjoint RITD$(3,t)$s. The existence result of a
pair of disjoint ITD$(3,t)$s and that of disjoint RITD$(3,t)$s are
given as follows.

\begin{Lemma}\label{itd}{\rm}
For any integer $t\geq 4$, there exists a pair of disjoint
ITD$(3,t)s.$ For any integer $t\geq 4$ and $t\neq 6,10$, there
exists a pair of disjoint RITD$(3,t)s.$
\end{Lemma}
\proof By \cite{fu}, for any integer $t\geq 4$, there exists a pair
of disjoint idempotent Latin squares of order $t$. Equivalently,
there is a pair of disjoint ITD$(3,t)s$.

By \cite{crc}, for any integer $t\geq 4$ and $t\neq 6, 10$, there
exist three mutually orthogonal Latin squares defined on and indexed
by $I_t$. By some permutations of rows and columns, we can form
three new mutually orthogonal Latin squares, say $L_1, L_2,L_3$, in
such a way that the main diagonal entries of $L_3$ are all 0's.
Accordingly, the main diagonal of $L_i$ ($i=1,2$) is a transversal.
By renaming the symbols of $L_1$ and $L_2$, we obtain two idempotent
Latin squares $L_1'$ and $L_2'$. Further $L_1', L_2'$ and $L_3$ are
still mutually orthogonal. Let $L_1'=(a_{ij})$, $L_2'=(b_{ij})$, and
$L_3=(c_{ij})$. For each $0\leq k\leq t-1$, let
$T_k=\{(i,j):c_{ij}=k\}$. Thus $T_0,T_1,\ldots, T_{t-1}$ form $t$
disjoint transversals of $L_1'$ and  $L_2'$, where $T_0$ consists of
the main diagonal positions. Then  we can construct a pair of
disjoint RITD$(3,t)s$ on $X= I_t\times I_3$ with group set ${\cal
G}=\{I_t\times \{i\}:i\in I_3\}$. For $0\leq k\leq t-1$, let
$P_1^k=\{\{(i,0),(j,1),(a_{ij},2)\}:(i,j)\in T_k\},$ and
$P_2^k=\{\{(i,0),(j,1),(b_{ij},2)\}:(i,j)\in T_k\}.$ It is readily
checked that  each $P_j^k$ ($0\leq k\leq t-1$, $j=1,2$) is a
parallel class of $X$ and $P_1^0=P_2^0$ is an idempotent parallel
class. Let ${\cal B}_1=\cup_{0\leq k\leq t-1}P_1^k$ and ${\cal
B}_2=\cup_{0\leq k\leq t-1}P_2^k.$ Observing that $a_{ij}\neq
b_{ij}$ if $i\neq j$, we obtain two disjoint RITD$(3,t)s$ $(X,{\cal
G},{\cal B}_1)$ and $(X,{\cal G},{\cal B}_2)$.\qed

%Some cases  have already been established. Here are some known
%results.
We next record some known results on disjoint 3-GDDs for later use.

\begin{Lemma}\label{intersect}

\begin{enumerate}
\item[{\rm (1)}]{\rm(\cite{butlerhoffman})} Let $u=0$, $g,t,u$  satisfy all
 the conditions of Lemma \ref{2.3}, and $(g,t,u)\neq(1,3,0)$. Then there exists a pair of disjoint $3$-GDDs of type $g^t$.
\item[{\rm (2)}]{\rm(\cite{lindner-rosa})} There exists a  pair of disjoint $3$-GDDs of type
$1^t3^1$, where $t\equiv 0,4$ {\rm(mod $6$)} and $t\geq  4$.
\end{enumerate}
\end{Lemma}
It is trivial that there  is a pair of disjoint $3$-GDDs of type
$g^t u^1$ if $gt=0$. And  Lemma \ref{intersect} solves the case
$u=g$ or  $u=0$. So we only need to consider the case $g,u$ all
positive, $u\neq g$, and $t\geq 3$. We call a triple $(g,t,u)$ of
positive integers with $u\neq g$ and $t\geq 3$ {\it admissible}
provided that the five conditions in Lemma \ref{2.3} all hold.

We shall utilize various methods to construct a pair of disjoint
$3$-GDDs of type $g^t u^1$ for any admissible triple $(g,t,u)$. And
we finally prove that  the necessary
 conditions for the existence  of a pair of $3$-GDDs of
type $g^tu^1$ are also sufficient. Our main result is:
\begin{Theorem}
\label{main}{\rm (Main Theorem)}  Let $g$, $t$, and $u$ be
nonnegative integers. There exists a a pair of disjoint $3$-GDDs of
type $g^tu^1$ if and only if the following conditions are all
satisfied:
\begin{enumerate}
\item[{\rm (1)}]  if $g>0$, then $t\geq3$ and $(g,t,u)\neq (1,3,0)$, or $t=2$ and $u=g$,
or $t=1$ and $u=0$, or $t=0;$

\item[{\rm (2)}] $u\leq g(t-1)$ or $gt=0;$

\item[{\rm (3)}] $g(t-1)+u \equiv 0$ $({\rm mod}$ $2)$ or $gt=0;$

\item[{\rm (4)}] $gt \equiv 0$ $({\rm mod}$ $2)$ or $u=0;$

\item[{\rm (5)}] $\frac{1}{2}g^2t(t-1)+gtu \equiv 0$ $({\rm mod}$ $3).$
\end{enumerate}

\end{Theorem}

\section{Recursive constructions}

In this section we  shall present several powerful recursive
constructions for disjoint $3$-GDDs.

The following construction is a variation of Wilson's Fundamental
Construction in \cite{wilson}.
\begin{Construction}
\label{WeightingConstruction} {\rm (Weighting Construction)} Suppose
that $(X, {\cal G}, {\cal A})$ is a $K$-GDD, and let $\omega :\,
X\longmapsto Z^{+}\cup \{ 0 \}$ be a weight function. For every
block $A\in {\cal A}$, suppose that there is a pair of disjoint
$3$-GDDs of type $\{\omega(x): \ x\in A\}$. Then there exists a pair
of disjoint $3$-GDDs of type $\{\sum_{x\in G}\omega(x): \ G\in {\cal
G} \}$.
\end{Construction}

\proof For every $x\in X$, let $S(x)$ be a set of $\omega(x)$
``copies'' of $x$. For any $Y\subseteq X$, let $S(Y)=\bigcup_{x\in
Y}S(x)$. For every block $A\in {\cal A}$, construct a pair of
disjoint $3$-GDDs $(S(A), \{S(x): x\in A\}, {\cal B}_A)$ and $(S(A),
\{S(x): x\in A\}, {\cal B}'_A\}$. Then it is readily checked that
there exists a pair of disjoint $3$-GDDs $(S(X)$, $\{S(G): G\in
{\cal G}\}, \cup_{A\in {\cal A}} {\cal B}_A)$ and $(S(X), \{S(G):
G\in {\cal G}\}, \cup_{A\in {\cal A}} {\cal B}'_A)$. \qed

We also employ ``Filling Construction" to break up the groups as
follows:

\begin{Construction}
\label{FillingGroups} {\rm (Filling Construction I)} Suppose that
there is a pair of disjoint $3$-GDDs of type $\{g_1, g_2, \ldots,$
$g_t\}$.  For each $1\leq i\leq t-1$, if $g_i\equiv 0$ {\rm(mod \
$s$)} and there is a pair of disjoint $3$-GDDs of type $s^{g_i/s
}u^1$. Then there exists a pair of disjoint $3$-GDDs of type
$s^{\sum_{i=1}^{t-1}g_i/s} (g_t+u)^1$.
\end{Construction}

%Further if $g_t\equiv 0$ {\rm(mod \ $s$)} and  a pair of disjoint
%$K$-GDDs of type $s^{g_t/s} u^1$ exists, then a pair of disjoint
%$3$-GDDs of type $s^{\sum_{i=1}^{t}g_i/s} u^1$ also exists.

\proof Let $(X,{\cal H},{\cal B}_1)$ and $(X,{\cal H},{\cal B}_2)$
be a pair of disjoint $3$-GDDs of type $\{g_1$, $g_2,$
$\ldots,g_t\}$. Let ${\cal H}=\{H_1, H_2, \ldots, H_t\}$ with
$|H_i|=g_i$ for $1\leq i\leq t$, and $Y$ be a set of cardinality $u$
such that $X\cap Y=\emptyset$.

For each $1\leq i\leq t-1$, we partition each $H_i$ into $g_i/s$
subsets $H_{ij}$, $1\leq j\leq g_i/s,$ such that $|H_{ij}|=s.$ By
assumption, there is a pair of $3$-GDDs on $H_i\bigcup Y$ with
$\{H_{ij}:1\leq j\leq g_i/s\}\cup\{Y\}$ as group set and ${\cal
A}_i^1$ and ${\cal A}_i^2$ as the disjoint block sets. Let ${\cal
G}=\{H_{ij}:1\leq i\leq t-1,1\leq j\leq g_i/s\}\cup\{H_t\cup Y\}$.
It is readily checked that $(X\bigcup Y, {\cal G}, (\cup_{i=1}^{t-1}
{\cal A}_i^1)\bigcup{\cal B}_1)$ and $(X\bigcup Y, {\cal G},
(\cup_{i=1}^{t-1} {\cal A}_i^2)\bigcup{\cal B}_2)$ are two disjoint
$3$-GDDs of type $s^{\sum_{i=1}^{t-1}g_i/s} (g_t+u)^1$.\qed

% Further if $g_t\equiv 0$ {\rm(mod \ $s$)}, then go on similarly to partition
%$H_t=\cup_{j=1}^{g_t/s} H_{tj}$ with $|H_{tj}|=s$ and construct on
%$H_t\cup \{Y\}$ a pair of $K$-GDDs of type $(s)^{g_t/s} u^1$ with
%disjoint block sets ${\cal A}_t^1$ and ${\cal A}_t^2.$ Thus we
%obtain a pair of  3-GDDs of type $s^{\sum_{i=1}^{t}g_i/s} u^1$ with
%disjoint block sets $(\cup_{i=1}^{t} {\cal A}_i^1)\bigcup{\cal B}_1$
%and $(\cup_{i=1}^{t} {\cal A}_i^2)\bigcup{\cal B}_2$.

\begin{Corollary} \label{2gtog} Let $t\geq 6$ be an even integer. If there exists a
pair of disjoint $3$-GDDs of type $(2g)^{t/2}u^1$, where $(g,
t/2)\neq (1,3)$, then so does  a pair of disjoint $3$-GDDs of type
$g^{t}(u+g)^1$.
\end{Corollary}

\proof It follows from Filling Construction I since  a pair of
disjoint $3$-GDDs of type $g^3$ exists by Lemma \ref{intersect}.
\qed

Sometimes we only fill in one long group and use the following
construction.
\begin{Construction}
\label{FillingOneGroup} {\rm (Filling Construction II)} Suppose that
there is a pair of disjoint $3$-GDDs of type $g^tu^1$ and $u=sg+x$.
If a pair of disjoint $3$-GDDs of type $g^sx^1$ also exists,  then
there exists a pair of disjoint $3$-GDDs of type $g^{s+t}x^1$.
\end{Construction}
\proof Let $(X,{\cal H}\cup\{G\},{\cal B}_1)$ and $(X,{\cal
H}\cup\{G\},{\cal B}_2)$ be a pair of disjoint $3$-GDDs of type $g^t
u^1$, where ${\cal H}=\{H_1,H_2,\ldots,H_t\}$ and $G=(\cup_{i=1}^{s}
G_i)\cup G_{s+1}$ with $|G_i|=g$ $(1\leq i\leq s$), $|G_{s+1}|=x$,
and $|H_j|=g$ $(1\leq j\leq t$). Construct on $G$ a pair of 3-GDDs
of type $g^s x^1$ with same group set ${\cal G}=\{G_i:1\leq i\leq
s+1\}$ and disjoint block sets ${\cal A}_1$ and ${\cal A}_2$. It is
immediately checked that $(X,{\cal G}\cup{\cal H}, {\cal A}_1\cup
{\cal B}_1)$ and $(X,{\cal G}\cup{\cal H}, {\cal A}_2\cup {\cal
B}_2)$ are two disjoint 3-GDDs of type $g^{s+t}x^1$. \qed

%\begin{Construction}
%\label{FillingSubdesigns} {\rm (Filling Construction)} Let $a$ be a
%nonnegative integer. Suppose that there exists a pair of disjoint
%$3$-GDDs and type $\{g_1, g_2, \ldots, g_s\}$. If there is a pair of
%disjoint $3$-GDDs of type $1^{g_i}a^1$ for each $1\leq i\leq s-1$,
%then there exists a pair of disjoint $3$-GDDs of type $1^{v-g_s}
%(g_s+a)$, where $v=\sum_{i=1}^s g_i$.
%\end{Construction}

What follows is a useful construction for generating 3-GDDs of type
$g^t u^1$ with $g$  relatively large.
\begin{Construction}
\label{glarge}   Suppose that there exists a $3$-GDD of type $\{g_1,
g_2, \ldots, g_s\}$. Let $t\geq 4$. If there is a pair of disjoint
$3$-GDDs of type ${g_i}^t u^1$ for each $1\leq i\leq s$, then there
exists a pair of disjoint $3$-GDDs of type $v^{t} u^1$, where
$v=\sum_{i=1}^s g_i$.
\end{Construction}
\proof Let $(X,{\cal G},{\cal B})$ be a 3-GDD of type $\{g_1, g_2,
\ldots, g_s\}$ and  $U$ be a set of cardinality $u$. We will
construct the desired designs on $(X\times I_t)\cup U$ with group
set ${\cal H}=\{X\times\{i\}:i\in I_t\}\cup\{U\}$.

For each block $B=\{x,y,z\}\in{\cal B}$,  there is a pair of
disjoint ITD$(3,t)$s by Lemma \ref{itd} on $B\times I_t$ with groups
$\{a\}\times I_t$, $a\in B$. Delete the idempotent parallel class to
form two disjoint block sets ${\cal A}_B^1$ and ${\cal A}_B^2$.

For each group $G\in {\cal G},$ place on $(G\times I_t)\cup U$ a
pair of disjoint 3-GDDs of type $|G|^t u^1$ with group set
$\{G\times\{i\}:i\in I_t\}\cup\{U\}$ and block sets ${\cal C}_G^1$
and ${\cal C}_G^2$.

Then we produce on  $(X\times I_t)\cup U$ a pair of disjoint 3-GDDs
of type $v^{t} u^1$ with  block sets $(\cup_{B\in {\cal B}}{\cal
A}_B^1)\cup(\cup_{G\in{\cal G}}{\cal C}_G^1)$ and $(\cup_{B\in {\cal
B}}{\cal A}_B^2)\cup(\cup_{G\in{\cal G}}{\cal C}_G^2).$\qed

\section{Direct constructions and preliminary results}

In this section we shall involve some methods of direct
construction. The ``method of differences" will be used to construct
some 3-GDDs of type $g^t u^1$, as is usually used in constructing
cyclic designs.  The cyclic partial Steiner triple systems also play
a crucial role  in constructing 3-GDDs.

The following result is simple but useful.

\begin{Lemma}\label{1-factor}
Suppose that there exists a pair of disjoint partial $3$-GDDs of
type $g^t u^1$ on $X$, where $U\subseteq X$ is the group of size
$u$, and $L_1$, $L_2$ are their leaves respectively. If the pairs of
the leave $L_j$ $(j=1,2)$ can be partitioned into $s$ disjoint
$1$-factors of $X\setminus U$, say, $F_1^j,F_2^j,\ldots,F_s^j$, such
that $F_i^1\cap F_i^2=\emptyset$ holds for each $1\leq i\leq s$,
then there exists a pair of disjoint $3$-GDDs of type $g^t(u+s)^1$.
\end{Lemma}

\proof Let $(X,{\cal G},{\cal B}_1)$ and $(X,{\cal G},{\cal B}_2)$
be  the assumed pair of disjoint partial $3$-GDDs of type $g^t u^1$
with $U$ as the group of size $u$. Define
$V=\{\infty_1,\infty_2,\ldots,\infty_s\}$, ${\cal
C}_j=\cup_{i=1}^s\{\{\infty_i,x,y\}:\{x,y\}\in F_i^j\}$, and ${\cal
H}=({\cal G}\setminus\{U\})\cup\{U\cup V\}$. Then $(X\cup V,{\cal
H},{\cal B}_1\cup{\cal C}_1)$ and $(X\cup V,{\cal H},{\cal
B}_1\cup{\cal C}_2)$ are two disjoint $3$-GDDs of type $g^t(u+s)^1$.
\qed

Each edge $\{a,b\}$ of a graph on vertices $Z_{v}$ is assigned to an
integer $d$ between 1 and $[v/2]$, called its {\it difference}, if
$|b-a|=d$ or $v-|b-a|=d$. A {\it difference triple} in $Z_v$ is a
set $\{a,b,c\}$ where $a+b\equiv c$ (mod $v$) or $a+b+c\equiv 0$
(mod $v$). A difference $d$ is called {\it good} in $Z_v$ if
$v/gcd(d,v)$ is even.

\begin{Lemma}\label{goodd}$({\rm \cite{stern}})$
Let $v$ be even and $D$ a subset of $[1,v/2]$. If $D$ contains a
good difference in $Z_v$, then the set of all unordered pairs of
$Z_v$ whose difference appears in $D$ can be partitioned into
1-factors.
\end{Lemma}

\begin{Lemma}\label{dif}
Let $(g,t,u)$ be an admissible triple with $u\geq 2$ and
$g(t-1)-u\equiv 0$ {\rm(mod\ $6$)}. Suppose that
$\{1,2,\ldots,gt/2\}\setminus\{t,2t,\ldots,[g/2]t\}=D_1\cup D_2$,
where $D_1$ can be partitioned into $(gt-g-u)/6$ difference triples
in $Z_{gt}$ and $gt/2\in D_2$ if $g$ is odd, or $D_2$ contains a
good difference in $Z_{gt}$ if $g$ is even, then there exists a pair
of disjoint $3$-GDDs of type $g^tu^1$.
\end{Lemma}

\proof Take $X=Z_{gt}\cup\{\infty_1,\infty_2,\ldots, \infty_{u}\}$
as the point set and ${\cal G}=\{\{j,t+j,2t+j,\ldots, (g-1)t+j\}:
0\leq j\leq t-1\}\cup\{\{\infty_1,\infty_2,\ldots, \infty_{u }\}\}$
as the group set. Suppose that $D_1$ can be partitioned into
difference triples $\{a_i,b_i,c_i\}$ in $Z_{gt}$ such that
$a_i+b_i\equiv c_i$ (mod $v$) or $a_i+b_i+c_i\equiv 0$ (mod $v$),
$1\leq i\leq (gt-g-u)/6$. Let
$$
{\cal A}_1=\cup_{1\leq i\leq (gt-g-u)/6}\{\{x,a_i+x,c_i+x\}:x\in
Z_{gt}\},
$$
and
$${\cal A}_2=\cup_{1\leq i\leq
(gt-g-u)/6}\{\{x,b_i+x,c_i+x\}:x\in Z_{gt}\}.
$$
Then $(Z_{gt}, {\cal A}_1)$ and $(Z_{gt}, {\cal A}_2)$ form two
disjoint partial 3-GDDs of type $g^t$. Their common leave $\cal L$
consists of all the pairs whose differences lie in $D_2$. By the
assumption, $D_2$ contains a good difference in $Z_{gt}$. By Lemma
\ref{goodd}, noting that $g$ and $u$ are both even or both odd,
$\cal L$ can be partitioned into $u$ $1$-factors, say, $F_1,
F_2,\ldots, F_{u}$. Let $F'_i=F_{i+1}$ for $i=1,2,\ldots,u$, where
the subscripts are modulo $u$. Since $u\geq 2$, $F_i\cap
F'_i=\emptyset$ $i=1,2,\ldots,u$. Hence, there exists a pair of
disjoint $3$-GDDs of type $g^tu^1$ by Lemma \ref{1-factor}. \qed

\begin{Corollary}\label{u=g(t-1)}
Let $u=g(t-1)$, where $g$ and $t$ are positive integers such that
$gt$ is even. Then there exists a pair of disjoint $3$-GDDs of type
$g^t u^1$.
\end{Corollary}

\proof The conclusion follows immediately by applying Lemma
\ref{dif} with $D_1=\emptyset$ and
$D_2=\{1,2,\ldots,gt/2\}\setminus\{t,2t,\ldots,[g/2]t\}$. \qed

A partial S$(2,3,v)$ is called {\it cyclic} if it has an
automorphism of order $v$. Usually, $Z_v$ is taken as the point set
of a cyclic design of order $v$ and the corresponding automorphism
is $i\rightarrow i+1$ (mod $v$). So the blocks of a partial
S$(2,3,v)$ can be partitioned into a number of orbits, each of which
can be represented by a {\it starter block}. An orbit is called {\it
full} if it consists of $v$ different blocks and called {\it short}
otherwise. In the proof of \cite[Lemma 3.2]{chr}, some cyclic
partial Steiner triple systems are constructed.

\begin{Lemma}\label{cps}{\rm (\cite{chr})}
For $k\geq 1$ and $1\leq s\leq 6,$  let $r'=7$ if $s=2$ and $k\equiv
2,3$ {\rm(mod\ $4$)}, or $r'=s-1$ otherwise. Then there is a cyclic
partial S$(2,3,6k+s)$ without short orbits whose leave is
$r$-regular, where
$r\equiv r'$ {\rm(mod\ $6$)}, $r'\leq r\leq 6k+s-1$. %For
%$k\geq 1$, there is a cyclic S$(2,3,6k+3)$.
Further if $r<6k+s-1$, then the cyclic partial S$(2,3,6k+s)$ has a
starter block containing a good difference.

\end{Lemma}

\begin{Lemma}\label{u>g}
Suppose that $(g,t,u)$ is an admissible triple with $u\geq 2$ and
$g(t-1)-u\equiv 0$ {\rm(mod\ $6$)}. Further suppose  $gt=6k+s$,
where $k\geq 1$ and $1\leq s\leq 6.$ Let  $r=7$ if $s=2$ and
$k\equiv 2,3$ {\rm(mod\ $4$)}, or $r=s-1$ otherwise. Whenever $u\geq
2g+r-2$ if $g$ is odd, or $u\geq 2g+r-5$ if $g$ is even, there
exists a pair of disjoint $3$-GDDs of type $g^tu^1$.
\end{Lemma}

\proof By Lemma \ref{cps}, there is a cyclic partial S$(2,3,gt)$
without short orbit whose leave is $r$-regular. Moreover, it has a
starter block containing a good difference. Let ${\cal F}$ be the
set of difference triples associated with the starter blocks of this
cyclic partial S$(2,3,gt)$. Let ${\cal F}_0$ be the set of
difference triples of ${\cal F}$, each of which contains at least a
multiple of $t$. Since $gt/2$ does not appear in a difference triple
of the cyclic partial S$(2,3,gt)$, we have $|{\cal F}_0|\leq
[(g-1)/2]$. Choose a subset ${\cal F'}$ such that ${\cal F}_0\subset
{\cal F'}\subset {\cal F}$ and $|{\cal F'}|=[(g-1)/2]$. Further for
even $g$ we can ensure that ${\cal F'}$ contains a difference triple
which have a good difference not being a multiple of $t$. This can
be done obviously if all the multiples of $t$ appear in less than
$(g-2)/2$ difference triples. Even if each difference triple of
${\cal F'}$ contains a multiple of $t$ as a difference, it can be
verified that the difference triple containing $t$ also contains a
good difference not being a multiple of $t$. Set $D_1=\cup_{ B\in
{\cal F}\setminus{\cal F'}} B$ and let $D_2$ be the set of
differences (between 1 and $gt/2$) neither appear in ${\cal
F}\setminus{\cal F'}$ nor are multiples of $t$. Since the cyclic
partial S$(2,3,gt)$ has no short orbit, we then have $D_1\cup D_2=
\{1,2,\ldots,gt/2\}\setminus\{t,2t,\ldots,[g/2]t\}$. Furthermore,
$|D_2|= g+(r-1)/2$ and $gt/2\in D_2$ if $g$ is odd, or $|D_2|=
g-2+(r-1)/2$ and $D_2$ contains a good difference in $Z_{gt}$ if $g$
is even. By Lemma \ref{dif}, there exists a pair of disjoint
$3$-GDDs of type $g^tu^1$, where $u=2g+r-2$ if $g$ is odd and $u=
2g+r-5$ if $g$ is even. For other cases of larger $u$ with
$g(t-1)-u\equiv 0$ {\rm(mod\ $6$)}, diverting more differences
produced by the difference triples in ${\cal F}\setminus {\cal F'}$
to $D_2$ works similarly. \qed

 Similar to Lemmas \ref{1-factor}, \ref{dif}, and
\ref{u>g}, we can obtain the result of disjoint partial $3$-GDDs of
type $g^{t} u^1$, whose leaves are same, forming a 1-factor of the
$t$ groups of size $g$. We record this in a remark.

{\Remark\label{remark}
 Suppose that $(g,t,u)$ is an admissible triple with   $u\neq 2$
 and $g(t-1)-u\equiv 0$ {\rm(mod\ $6$)}.  Further suppose
$gt=6k+s$, where $k\geq 1$ and $1\leq s\leq 6.$ Let  $r=7$ if $s=2$
and $k\equiv 2,3$ {\rm(mod\ $4$)}, or $r=s-1$ otherwise.  Whenever
$u\geq 2g+r-2$ if $g$ is odd, or $u\geq 2g+r-5$ if $g$ is even,
there exists a pair of disjoint partial $3$-GDDs of type
$g^t(u-1)^1$, whose leaves are same, forming a $1$-factor of the $t$
groups of size $g$. }

Next we consider two small cases $g=1$ and $g=2$.

\begin{Lemma}
\label{g11}$({\rm \cite{chu}})$ There exists a pair of disjoint
$3$-GDDs of type $1^t u^1$ whenever $u\equiv 1, 3$ {\rm(mod $6$)},
$u+t\equiv 1, 3$ {\rm(mod $6$)} and $7\leq u\leq t-1$.
\end{Lemma}

%Let $X$ be a $u$-set and $Y$ be a $v$-subset of $X$, where $u,
%v\equiv 1, 3$ \rm(mod $6$), $u\geq 2v+1$ and $v\geq 7$. There exist
%a pair of disjoint STS$(u)$s $(X,{\cal B}_1)$ and $(X,{\cal B}_2)$
%and a pair of disjoint STS$(v)$s $(Y,{\cal C}_1)$ and $(Y,{\cal
%C}_2)$ such that ${\cal C}_1\subset {\cal B}_1$ and ${\cal
%C}_2\subset {\cal B}_2$.

%Then we can completely solve the case $g=1.$

\begin{Lemma}\label{g1}%$({\rm \cite{}})$
The Main Theorem holds for any admissible triple $(1,t,u)$.
\end{Lemma}

 \proof Since $(1,t,u)$ is an admissible triple, $u$ must be odd and $u\geq 3$. We distinguish the possibility of $u$ to show the conclusion.

 First if $u=3$, then  $t\equiv 0,4 $ (mod 6)
and $t\geq 4$. A pair of disjoint 3-GDDs of type $1^{t}3^1$ exists
by Lemma \ref{intersect}.

Next if $u\equiv 1, 3$ (mod $6$) and $u\geq 7$, then $u+t\equiv 1,
3$ (mod $6$) and $u\leq t-1$.   By Lemma \ref{g11}, there exists a
pair of disjoint $3$-GDDs of type $1^t u^1$.

Finally we  treat $u\equiv 5 $ (mod 6). Then $t\equiv 0 $ (mod 6)
and $u\leq t-1$. Corollary \ref{u=g(t-1)} solves the case $t=6$ and
$u=5$. For $t\geq 12$, a pair of disjoint $3$-GDDs of type $1^{t}
u^1$ is obtained by taking $g=1$ and $r=5$ in Lemma \ref{u>g}.
  \qed

%For $t=4,6,10,12$, a pair of disjoint 3-GDDs of type $1^{t}3^1$
%exists by Lemma \ref{intersect}. For $t\geq 16$, since there are a
%pair of disjoint $3$-GDDs of type $1^{t-6} 9^1$ and a pair of
%disjoint $3$-GDDs of type $1^6 3^1$ by Lemma \ref{intersect}, we can
%obtain a pair of disjoint $3$-GDDs of type $1^{t} 3^1$ by Filling
%Construction II.

% So we can treat the case $g=2$.

\begin{Lemma}\label{g21}%$({\rm \cite{}})$
The Main Theorem holds for any admissible triple $(2,t,u)$ with
 $t\equiv 1,2$ {\rm(mod $3$)}.
\end{Lemma}

\proof  Since $(2,t,u)$ is an admissible triple, $t\equiv 1$
{\rm(mod $3$)} requires $u\equiv 0$ {\rm(mod $6$)} ($u\geq 6$),
$t\equiv 2$ {\rm(mod $3$)} demands $u\equiv 2$ {\rm(mod $6$)}
($u\geq 8$), and $(1,2t,u+1)$ is also an admissible triple
satisfying the equality $1\cdot(2t-1)-(u+1)\equiv 0$ (mod 6). Let
$2t=6k+s$ and $k,s,r$ be taken as  in Remark \ref{remark}.
 As  $u+1\geq 7\geq r=2\cdot
1+r-2$,  there is a pair of partial $3$-GDDs of type $1^{2t} u^1$
with $U$ as the long group, whose leaves are same, forming a
$1$-factor of the $2t$ groups of size $1$. Take this 1-factor
together with $U$ as new groups, we obtain a pair of disjoint
$3$-GDDs of type $2^t u^1$. \qed

The complete solution for the case $g=2$ is left to Section 5.

\section{The  case $t\equiv 3$ $({\rm mod}$ $6)$}

%\section{Special pairs of resolvable $\{2,3\}$-GDDs }

A useful auxiliary design to construct 3-GDDs is  resolvable
$\{2,3\}$-GDD with 3 groups of even size, whose existence is
investigated in \cite{rees}. We shall show in this section that two
such GDDs with some restrictions also exist. Related results will be
employed to solve the case $t\equiv 3$ $({\rm mod}$ $6)$ of the Main
Theorem.

\begin{Lemma}\label{23gdd}
Let $g$ and $u$ be even, $0\leq u\leq 2g$, $(g,u)\neq(2,0)$ or
$(6,0)$. Then there is a pair of $\{2,3\}$-GDD of type $g^3$ with
same groups and different block sets ${\cal B}^1$ and ${\cal B}^2$
satisfying all of the following conditions{\rm :}

\begin{enumerate}
\item[{\rm (1)}] Both ${\cal B}^1$ and ${\cal B}^2$ can be resolved
into $u$ parallel classes containing only blocks of size $2$ and
$g-u/2$ parallel classes containing only blocks of size $3;$

\item[{\rm (2)}]  ${\cal B}^1$ and ${\cal B}^2$ have no block of size
$3$ in common{\rm ;}

\item[{\rm (3)}] The $u$ parallel classes containing only blocks
 of size $2$ of  ${\cal B}^j$ $(j=1,2)$ can be arranged in sequence
  $P_1^j, P_2^j,\ldots,P_u^j,$ in such a way that $P_i^1\cap P_i^2=\emptyset$ for each $1\leq i\leq u$.
\end{enumerate}

\end{Lemma}

\proof We follow the idea of Rees in \cite{rees}. Let $X=Z_g\times
I_3$ be the point set and ${\cal G}=\{Z_g\times\{i\}$: $i\in I_3\}$
be the group set.

First we handle the case $u=0$. Obviously when $g\neq 2,6$, there
exists a resolvable 3-GDD $(X, {\cal G}, {\cal C})$ of type $g^3$.
Set ${\cal
C}'=\{\{(x,0),(y,1),(z+1,2)\}:\{(x,0),(y,1),(z,2)\}\in{\cal C}\}$.
Then $(X, {\cal G}, {\cal C}')$ is a resolvable 3-GDD disjoint with
$(X, {\cal G}, {\cal C})$.

Next consider $u\geq 2$.
 Let ${\cal B} $ be the union of following $g+1$ parallel
classes of $X$:
\begin{eqnarray*}
& & S_i =\{\{(x,0),(x+i,1),(x+2i,2)\}:x\in Z_g\}, \ 0\leq i\leq g/2-1,\\
& & S_i =\{\{(x,0),(x+i,1),(x+2i+1,2)\}:x\in Z_g\}, \ g/2\leq i\leq g-2,\\
& & M_1 =\{\{(x,0),(x-1,1)\},\{(x+g/2,0),(x+g/2-1,2)\},\\
&&\ \ \ \ \ \ \ \ \ \ \{(x+g/2-1,1),(x-1,2)\}:
0\leq x\leq g/2-1\},\\
& & M_2 =\{\{(x,0),(x-1,1)\},\{(x+g/2,0),(x+g/2-1,2)\},\\
&&\ \ \ \ \ \ \ \ \ \ \{(x+g/2-1,1),(x-1,2)\}:g/2\leq x\leq g-1\}.
\end{eqnarray*}
Then $(X,{\cal G},{\cal B} )$ is a resolvable $\{2,3\}$-GDD
 with two
parallel classes of blocks of size 2.

To generate more parallel classes, some transformations from
parallel classes of triples to those of pairs are made.

 (A)  The pairs produced by
$S_{g/2-1} $ and $M_1 $ can be divided into three parallel classes
$P_{1l} $, $1\leq l\leq 3$, described below.  Let\begin{eqnarray*}
& & M_{11} =\{\{(x,0),(x-1,1)\}:0\leq x\leq g/2-1\}, \\
& & M_{12} =\{\{(x,0),(x-1,2)\}:g/2\leq x\leq g-1 {\rm\ and\ }x{\rm \ is\ even}\}\\
&&\ \ \ \ \ \ \ \ \ \ \cup\{\{(x+g/2-1,1),(x-1,2)\}:
0\leq x\leq g/2-1  {\rm\ and\ }x{\rm \ is\ even}\},\\
& & M_{13} =(M_1 \setminus  M_{11} )\setminus  M_{12} .
\end{eqnarray*}
For each block $B$ of $S_{g/2-1} $ and $1\leq l\leq 3$, let $h_l^1
(B)$ be the unique intersection of $B$ and $M_{1l} $ and let
$$ P_{1l} =M_{1l} \cup(\cup\{B\setminus
\{h_l^1 (B)\}:B\in S_{g/2-1}\}).$$

Note: By  replacing $M_1$  with $M_2$ and ``$x$ is even" with ``$x$
is odd" and interchanging the range $0\leq x\leq g/2-1$ and $g/2\leq
x\leq g-1$ in $M_{1l}$, the pairs produced by $S_{g/2-1} $ and $M_2
$ can also be divided into three parallel classes, which we denote
by $P_{2l} $, $1\leq l\leq 3$.

 (B) For $0\leq i\leq g/2-2$, all the pairs
produced by the two classes $S_i $ and $S_{g/2+i} $ can be divided
into four parallel classes $E_{ik} $, $1\leq k\leq 4$, as follows:
\begin{eqnarray*}
& & E_{i1} =\{\{(2x,0),(2x+i,1)\},\{(2x+1,0),(2x+2i+2,2)\},\\
& & \ \ \ \ \ \ \ \ \ \ \{(2x+i+1,1),(2x+2i+1,2)\}:0\leq x\leq g/2-1\}, \\
& & E_{i2} =\{\{(2x+1,0),(2x+g/2+i+1,1)\},\{(2x,0),(2x+2i,2)\},\\
& & \ \ \ \ \ \ \ \ \ \ \{(2x+g/2+i,1),(2x+2i+1,2)\}:0\leq x\leq
g/2-1\}.
\end{eqnarray*}
Setting $E_{i,k+2}=\{\{(x+1,s),(y+1,t)\}:\{(x,s),(y,t)\}\in
E_{ik}\}$ for $k=1,2$ yields another two parallel classes $E_{i3}$
and $E_{i4}$.

  Let $\phi$ be a bijection on $Z_g\times I_3$ such
that $\phi((x,0))=(x,0),$ $\phi((x,1))=(x,1),$ and
$\phi((x,2))=(x+1,2).$ For a subset ${\cal A}$ of ${\cal B},$ define
$\phi({\cal A})=\{\{\phi(a),\phi(b),\phi(c)\}:\{a,b,c\}\in {\cal
A}\}$.

 If $u/2$
is odd, then in ${\cal B} $  by replacing $S_i $ and $S_{g/2+i} $
with $E_{ik} $ (only if $u\geq 6$) for $0\leq i\leq (u-6)/4$, $1\leq
k\leq 4$, we obtain a  resolvable $\{2,3\}$-GDD $(X,{\cal G},{\cal
B}^1)$ with exactly $u$ parallel classes of pairs. ${\cal P}_1=\{M_l
:l=1,2\}\cup\{E_{ik}:0\leq i\leq (u-6)/4, 1\leq k\leq 4\}$ is the
collection of the $u$ parallel classes of pairs. And ${\cal
P}_2=\{S_i :(u-2)/4\leq i\leq g/2-1,$ or $(u-2)/4+g/2\leq i\leq
g-2\}$  is the collection of the parallel classes of triples. Let
${\cal P}={\cal P}_1\cup{\cal P}_2 $ and ${\cal B}^2=\phi({\cal
B}^1)$. Apparently, $(X,{\cal G},{\cal B}^2)$ is a resolvable
$\{2,3\}$-GDD with a collection of parallel classes $\{\phi(P):
P\in{\cal P}\}$. Besides, one can check that $\phi(M_1)\cap
M_2=\emptyset $, $\phi(M_2)\cap M_1=\emptyset $, $\phi(E_{ik})\cap
E_{i,k+2}=\emptyset$ $(0\leq i\leq (u-6)/4$, $k,k+2$ is modulo 4),
and $\phi(Q)\cap R=\emptyset$ for any $Q,R\in {\cal P}_2.$ So we
prove the lemma for $u/2$ odd.

Otherwise, $u/2$ is even. Then  in ${\cal B} $  by replacing $S_i $
and $S_{g/2+i} $ with $E_{ik} $ (only if $u\geq 8$) for $0\leq i\leq
(u-8)/4$, $1\leq k\leq 4$, and replacing $S_{g/2-1} $ and $M_1 $
with $P_{1l} $, $1\leq l\leq 3$, we obtain a   resolvable
$\{2,3\}$-GDD $(X,{\cal G},{\cal B}^1)$ with exactly $u$ parallel
classes of pairs. ${\cal P}_1=\{E_{ik}:0\leq i\leq (u-8)/4,1\leq
k\leq 4\}\cup\{M_2\}\cup\{P_{1l}:l=1,2,3\}$ contains the $u$
parallel classes of pairs. And ${\cal P}_2=\{S_i:(u-4)/4\leq i\leq
g/2-2,$ or $(u-4)/4+g/2\leq i\leq g-2\}$ contains all the parallel
classes of triples. If we employ the same replacement except taking
$M_2$ instead of $M_1$, then another resolvable $\{2,3\}$-GDD
$(X,{\cal G},{\cal B}')$ is obtained. The collection of parallel
classes are ${\cal P}'=(({\cal P}_1\cup{\cal P}_2)\setminus\{ M_2,
P_{11}, P_{12},P_{13}\})\cup\{M_1\}\cup\{P_{2l}:l=1,2,3\}$.  Let
${\cal B }^2=\phi({\cal B}').$  Then $(X,{\cal G},{\cal B }^2)$ is a
resolvable $\{2,3\}$-GDD of type $g^3$ with a collection of parallel
classes $\{\phi(P): P\in{\cal P}'\}$. Further, $ {\cal B}^1  $ and
${\cal B}^2$ satisfy  the three conditions required by the lemma,
where $\phi(E_{ik})\cap E_{i,k+2}=\emptyset $ ($0\leq i\leq
(u-8)/4,$ $k,k+2$ is modulo 4), $\phi(M_1)\cap M_2=\emptyset,$ and
$\phi(P_{2l})\cap P_{1l}=\emptyset $ $(l=1,2,3)$, $\phi(Q)\cap
R=\emptyset$ for any $Q,R\in {\cal P}_2.$ This completes the proof.
\qed

\begin{Corollary}\label{t3}
The Main Theorem holds for any admissible triple $(g,t,u)$ with
$t\equiv 3$ {\rm(mod\ $6$)}.
\end{Corollary}
\proof $(g,t,u)$ is admissible and $t\equiv 3$ {\rm(mod\ $6$)}, so
$g\equiv 0$ {\rm(mod\ $2$)}, $u\equiv 0$ {\rm(mod\ $2$)}, and $2\leq
u\leq g(t-1)$.

 We first treat $t=3.$ Suppose that
$(X,{\cal G},{\cal A}_1\cup{\cal B}_1)$ and $(X,{\cal G},{\cal
A}_2\cup{\cal B}_2)$ are
 two $\{2,3\}$-GDD of type $g^3$ satisfying all the three
conditions in Lemma \ref{23gdd}, where ${\cal A}_i$ $(i=1, 2)$
consists of $u$ parallel classes of pairs, say, $F_1^i,
F_2^i,\ldots,F_{u}^i$, and ${\cal B}_i$ $(i=1, 2)$ consists of
parallel classes of triples. Further $F_j^1\cap F_j^2=\emptyset$ for
$1\leq j\leq u$ and ${\cal B}_1\cap{\cal B}_2=\emptyset$. By Lemma
\ref{1-factor}, there is a pair of disjoint 3-GDDs of type $g^3u^1$.

 Next let $t=6n+3$ where $n\geq 1$. There is a
 KTS$(t)$ on a $t$-set $Y$ having $3n+1$ parallel classes
 $P_1,P_2,\ldots,P_{3n+1}$. Since $u\equiv 0$ {\rm(mod\ $2$)} and $u\leq g(t-1)$,
  we can take even integers $u_j$, $j=1, 2,\ldots,3n+1,$ such that
  $0\leq u_j\leq 2g$ and $u=\sum_{j=1}^{ 3n+1}u_j$. Let $U_j=\{\infty_1^j,\infty_2^j,\ldots,
 \infty_{u_j}^j\}$ and $U=\cup_{j=1}^{ 3n+1}U_j$. For every block $B=\{x,y,z\}$ of each parallel class $P_j$, $1\leq j\leq
 3n+1$, construct on $(B\times I_g)\cup U_j$ a pair of disjoint 3-GDDs of type $g^3{u_j}^1$
 with group set $\{\{x\}\times I_g:x\in B\}\cup\{U_j\}$ and block
 sets ${\cal C}_B^1$ and ${\cal C}_B^2$. Set $Z=(Y\times I_g)\cup U$, ${\cal G}=\{\{x\}\times I_g:x\in Y\}\cup\{U\}$
  and ${\cal C}^i=\bigcup_{B\in P_j,1\leq j\leq 3n+1}{\cal C}_B^i$ for $i=1, 2$. It is immediate
 that $(Z,{\cal G},{\cal C}^1)$ and $(Z,{\cal G},{\cal C}^2)$ are
 two disjoint 3-GDDs of type $g^tu^1$.
 \qed

\begin{Lemma}\label{23gdd2}
Let $g$ and $u$ be even, $2\leq u\leq 2g-2$. Then there is a pair of
$\{2,3\}$-GDD of type $g^3$ with same groups and different block
sets ${\cal B}^1$ and ${\cal B}^2$ satisfying all of the following
conditions{\rm :}

\begin{enumerate}
\item[{\rm (1)}] Both ${\cal B}^1$ and ${\cal B}^2$ can be resolved
into $u$ parallel classes containing only blocks of size $2$ and
$g-u/2$ parallel classes containing only blocks of size $3;$

\item[{\rm (2)}]  ${\cal B}^1$ and ${\cal B}^2$ have a common parallel class of size $3$ but have  no other triple in common{\rm ;}

\item[{\rm (3)}] The $u$ parallel classes containing only blocks
 of size $2$ of  ${\cal B}^j$ $(j=1,2)$ can be arranged in sequence
  $P_1^j, P_2^j,\ldots,P_u^j,$ in such a way that $P_i^1\cap P_i^2=\emptyset$ for each $1\leq i\leq u$.
\end{enumerate}
\end{Lemma}
\proof The proof is similar to that of Lemma \ref{23gdd}. First we
have a resolvable $\{2,3\}$-GDD $(X,{\cal G},{\cal B})$ of type
$g^3$ with $M_1$ and $M_2$ as  the parallel classes of pairs, and
$S_i$, $0\leq i\leq g-2,$ as the parallel classes of triples. The
conclusion holds clearly for the case $(g,u)=(2,2)$, so we assume
that $g\geq 4$. We will use transformation of kind (B) and another
three kinds to treat the parallel classes.

(C)  The pairs produced by $S_{0} $ and $M_1 $ can be divided into
three parallel classes $P_{0l} $, $1\leq l\leq 3$.
Let\begin{eqnarray*}
& & M_{01} =\{\{(x+g/2-1,1),(x-1,2)\}:0\leq x\leq g/2-1\}, \\
& & M_{02} =\{\{(x,0),(x-1,2)\}:g/2\leq x\leq g-1 {\rm\ and\ }x{\rm \ is\ even}\}\\
&&\ \ \ \ \ \ \ \ \ \ \cup\{\{(x,0),(x-1,1)\}:
0\leq x\leq g/2-1  {\rm\ and\ }x{\rm \ is\ even}\},\\
& & M_{03} =(M_1 \setminus  M_{01} )\setminus  M_{02} .
\end{eqnarray*}
For each block $B$ of $S_{0} $ and $1\leq l\leq 3$, let $h_l^0 (B)$
be the unique intersection of $B$ and $M_{0l} $ and let
$$ P_{0l} =M_{0l} \cup(\cup\{B\setminus
\{h_l^0 (B)\}:B\in S_{0}\}).$$

(D)  The pairs produced by the two classes $S_{0} $ and $S_{g-2} $
can be divided into four parallel classes $F_{k} $, $1\leq k\leq 4$,
as follows:
\begin{eqnarray*}
& & F_{1} =\{\{(2x+1,0),(2x-1,1)\},\{(2x,0),(2x,2)\},\\
& & \ \ \ \ \ \ \ \ \ \ \{(2x,1),(2x-1,2)\}:0\leq x\leq g/2-1\}, \\
& & F_{2}
=\{\{(2x,0),(2x,1)\},\{(2x+1,0),(2x-2,2)\},\\
& & \ \ \ \ \ \ \ \ \ \ \{(2x+1,1),(2x+1,2)\}:0\leq x\leq g/2-1\}.
\end{eqnarray*}
Setting $F_{k+2}=\{\{(x+1,s),(y+1,t)\}:\{(x,s),(y,t)\}\in F_{k}\}$
for $k=1,2$ yields another two parallel classes $F_{3}$ and $F_{4}$.

(E)  The pairs produced by the two classes $S_{g/2-2} $ and
$S_{g/2-1} $ can be divided into four parallel classes $H_{k} $,
$1\leq k\leq 4$, as follows:
\begin{eqnarray*}
& & H_{1} =\{\{(2x+1,0),(2x+g/2-1,1)\},\{(2x,0),(2x-4,2)\},\\
& & \ \ \ \ \ \ \ \ \ \ \{(2x+g/2,1),(2x-1,2)\}:0\leq x\leq g/2-1\}, \\
& & H_{2} =\{\{(2x,0),(2x+g/2-1,1)\},\{(2x+1,0),(2x-1,2)\},\\
& & \ \ \ \ \ \ \ \ \ \ \{(2x+g/2-2,1),(2x-4,2)\}:0\leq x\leq
g/2-1\}.
\end{eqnarray*}
Setting $H_{k+2}=\{\{(x+1,s),(y+1,t)\}:\{(x,s),(y,t)\}\in H_{k}\}$
for $k=1,2$ yields another two parallel classes $H_{3}$ and $H_{4}$.

 \indent  Let $\phi$ be a bijection on $Z_g\times I_3$
such that $\phi((x,0))=(x,0),$ $\phi((x,1))=(x+1,1),$ and
$\phi((x,2))=(x+3,2).$ For a subset ${\cal A}$ of ${\cal B}$ define
$\phi({\cal A})=\{\{\phi(a),\phi(b),\phi(c)\}:\{a,b,c\}\in {\cal
A}\}$. Evidently, $\phi(S_{g/2-1})=S_{g/2},$ which we will use as
the common parallel class required by the lemma.

First let $u/2$ be odd. If more parallel classes of pairs are
required, then replace step by step in ${\cal B} $ each pair $S_0 $
and $S_{g-2} $ with $F_{k} $, $S_{g/2-2} $ and $S_{g/2-1} $ with
$H_{k}$, $S_i $ and $S_{g/2+i} $ with $E_{ik} $ ($1\leq i\leq
(u-10)/4$, $1\leq k\leq 4$). Thus we obtain a resolvable
$\{2,3\}$-GDD $(X,{\cal G},{\cal B}^1)$ with a collection of
parallel classes
 ${\cal
P}={\cal P}_1\cup{\cal P}_2$, where ${\cal
P}_1=\{M_i:i=1,2\}\cup\{F_{k}:1\leq k\leq 4\}\cup\{H_{k}:1\leq k\leq
4\}\cup\{E_{ik}: 1\leq i\leq (u-10)/4, 1\leq k\leq 4\}$, ${\cal
P}_2=\{S_i: i=g/2,$ or $(u-6)/4\leq i\leq g/2-3,$ or $
(u-6)/4+g/2\leq i\leq g-3\}$ (observe that $S_{g/2}\in{\cal P}$).
Similarly, replace in ${\cal B} $ each pair $S_0 $ and $S_{g/2} $
with $E_{0,k} $, $S_{g/2-2} $ and $S_{g-2} $ with $E_{g/2-2,k}$. And
we still replace $S_i $ and $S_{g/2+i} $
 with $E_{ik} $ ($1\leq i\leq (u-10)/4$, $1\leq k\leq 4$), then form
 another resolvable $\{2,3\}$-GDD $(X,{\cal G},{\cal B}')$ with a
collection of parallel classes ${\cal P}'={\cal P}'_1\cup{\cal
P}'_2$, where ${\cal P}'_1=\{M_i:i=1,2\}\cup\{E_{ik}: 0\leq i\leq
(u-10)/4,$ or $i=g/2-2, 1\leq k\leq 4\}$, ${\cal P}'_2=\{S_i:
(u-6)/4\leq i\leq g/2-3,$ or $i=g/2-1$, or $ (u-6)/4+g/2\leq i\leq
g-3\}$. Let ${\cal B }^2=\phi({\cal B}').$ Obviously, $(X,{\cal
G},{\cal B }^2)$ is a resolvable $\{2,3\}$-GDD of type $g^3$ with a
collection of parallel classes $\{\phi(P): P\in{\cal P}'\}$
containing $\phi(S_{g/2-1})$. Besides, one can check that
$\phi(P)\cap P=\emptyset$ for any $P\in {\cal
P}'_1\setminus\{E_{0k},E_{g/2-2,k}:1\leq k\leq 4\}$,
$\phi(E_{0k})\cap F_k=\emptyset$, $\phi(E_{g/2-2,k})\cap
H_k=\emptyset$ (a slight difference when $g/2$ is odd:
$\phi(E_{g/2-2,2})\cap H_4=\phi(E_{g/2-2,4})\cap H_2=\emptyset$),
and $\phi(Q)\cap R=\emptyset$ for any $Q, R\in {\cal P}'_2$ except
$\phi(S_{g/2-1})=S_{g/2}.$

Finally let $u/2$ be even. For $1\leq i\leq (u-4)/4$, $1\leq k\leq
4,$ replace in ${\cal B} $ each pair $S_i $ and $S_{g/2+i} $ with
$E_{ik} $, and replace $S_{0} $ and $M_1 $ with $P_{0l} $, $1\leq
l\leq 3$. Thus we obtain a   resolvable $\{2,3\}$-GDD $(X,{\cal
G},{\cal B}^1)$ with  a collection of parallel classes ${\cal
P}={\cal P}_1\cup{\cal P}_2$, where ${\cal P}_1=\{E_{ik}: 1\leq
i\leq (u-4)/4, 1\leq k\leq 4\} \cup\{P_{0l} :l=1,2,3\}\cup\{M_2\}$,
${\cal P}_2=\{S_i: u/4\leq i\leq g/2,$ or $u/4+g/2\leq i\leq g-2\}$
(note that both $S_{g/2-1}$ and $S_{g/2}$ belong to ${\cal P}$).
Similarly let ${\cal B }^2=\phi({\cal B}^1).$ Then $(X,{\cal
G},{\cal B }^2)$ is a resolvable $\{2,3\}$-GDD  of type $g^3$ with a
collection of parallel classes $\{\phi(P): P\in{\cal P}\}$, which
also satisfy all the conditions required by the lemma.
 \qed

\begin{Corollary}\label{t3i}
Let $g$ and  $u$ be even integers such that $0\leq u\leq 2g-2$  and
$(g,u)\neq (2,0)$. Then there exists a pair of $3$-GDDs of type
$g^3u^1$ with exactly $g$ blocks in common and these  $g$ blocks
form a parallel class of the union of the three groups of size $g$.
\end{Corollary}
\proof There is a pair of disjoint ITD$(3,g)$s for $g\geq 4$ by
Lemma \ref{itd}, so the conclusion holds if $u=0$. If $2\leq u\leq
2g-2$, there is a pair of $\{2,3\}$-GDDs meeting the conditions in
Lemma \ref{23gdd2}. Analogous to the proof for $t=3$ in Corollary
\ref{t3}, the conclusion follows. \qed

\section{The  case $g\equiv 0$ $({\rm mod}$ $3)$}

In this section, we mainly examine the existence of a pair of
disjoint 3-GDDs of type $g^tu^1$ for $g\equiv 0$ ${\rm (mod}$ $3)$.
We adopt a similar procedure as in Section 2 of \cite{chr}, so we
list some results on $K$-GDDs derived therein.

\begin{Lemma}
\label{kgdd}$({\rm
\cite{assafhartman,chr,raywilson,rees,reesstinson}})$
\begin{enumerate}
\item[{\rm (1)}] For odd integer $t\geq 3$, there is a $4$-GDD of
type $3^t (\frac{3(t-1)}{2})^1.$

\item[{\rm (2)}] For even integer $t\geq 6$, there is a $\{4,7\}$-GDD of type $3^t
(\frac{3(t-2)}{2})^1,$ in which  precisely one point of the long
group belongs to blocks of size $7.$ Further   this point does not
belong to any block of size $4$ if $t\geq 8.$

\item[{\rm (3)}] There is a $4$-GDD of type $3^5.$

\item[{\rm (4)}]  For $(t,m,k)=(4,6,3),(6,8,1)$, there is a
$\{3,4\}$-GDD of type $3^t m^1$, in which precisely $k$ points of
the long group belong to the blocks of size $3.$
\end{enumerate}
\end{Lemma}

The following three lemmas  are all presented by utilizing the
Weighting Construction. So we only point out the initial $K$-GDDs
(all coming from Lemma \ref{kgdd}), the weight function, and the
input designs in the proof.

\begin{Lemma}\label{g6x1}
The Main Theorem holds for any admissible triple $(g,t,u)$ with
$g\equiv 0$ {\rm (mod $6$)} and $t\equiv 1$ {\rm (mod $2$)}.
 \end{Lemma}
 \proof Let $g=6x$ where $x\geq 1$. Start from  a 4-GDD of type
 $3^t(\frac{3(t-1)}{2})^1$ with a
 long group $Y=\{y_1,y_2,\ldots,y_{3(t-1)/2}\}$  Then give even weight $w_i$ between 0 and $4x$ to each
 point $y_i$
 of $Y$ such that $u=\sum_{i=1}^{3(t-1)/2}w_i$. Next give weight $2x$ to any other point. By Lemma \ref{intersect} and Corollary \ref{t3}, for even $0
 \leq w\leq 4x$ there is a pair of disjoint 3-GDDs of type
 $(2x)^3w^1$. So the conclusion follows by the Weighting Construction.
 \qed

\begin{Lemma}\label{g6x2}
The Main Theorem holds for any admissible triple $(g,t,u)$ with
$g\equiv 0$ {\rm (mod $6$)}, $t\equiv 0$ {\rm (mod $2$)}, and $t\geq
8$.
 \end{Lemma}
 \proof  Let $g=6x$ where $x\geq 1$. Start from a  $\{4,7\}$-GDD of type $3^{t}
 (\frac{3(t-2)}{2})^1$ with a long
 group $Y=\{y_1,y_2,\ldots,y_{3(t-2)/ 2}\}$, where  only one point $y_1$ of $Y$ belongs to the block of size
7,
  and $y_1$  does not belong to any block of size $4$. We  give $y_1$
weight $w_1=0$ or $10x$, give each $y_i\in Y$ with $i\geq 2$ even
weight $w_i,$
 $0 \leq w_i\leq 4x,$ such that $u=\sum_{i=1}^{3(t-2)/ 2}w_i$, and give each point not in $Y$ weight $2x$.  Since  two disjoint 3-GDDs of type
$(2x)^3w^1$ $(w$ even, $0\leq w\leq 4x)$,   or $(2x)^6 v^1$
($v=0,10x$) exist by  Lemma \ref{intersect}, Corollaries
\ref{u=g(t-1)} and \ref{t3}, a pair of disjoint $3$-GDDs of type
$g^tu^1$ is obtained.
 \qed

\begin{Lemma}\label{t46}
The Main Theorem holds for any admissible triple $(g,t,u)$ with
$g\equiv 0$ {\rm (mod $6$)} and $t=4,6$.
\end{Lemma}

\proof  Let $g=6x$ where $x\geq 1$. Set $(m,k)=(6,3)$ if $t=4$ and
$(m,k)=(8,1)$ if $t=6$.

First we handle even $u$ with $2kx\leq u\leq g(t-1)$. Start from a
$\{3,4\}$-GDD of type $3^t m^1$ with a long group
$Y=\{y_1,y_2,\ldots,y_m\}$, in which  precisely $k$ points
$y_1,y_2,\ldots,y_k$ belong to the blocks of size $3.$ Give each
$y_i$ with $1\leq i\leq k$ weight $2x$ and each $y_i$ with $k+1\leq
i\leq m$ even weight $w_i,0\leq w_i\leq 4x$ such that
$u=2kx+\sum_{i=k+1}^m w_i$. Then weight $2x$ to every point not in
$Y$. Since a pair of disjoint 3-GDDs of type $(2x)^3 w^1$ ($w$ even,
$0\leq w\leq 4x)$ exists by Lemma
 \ref{intersect} and Corollary \ref{t3}, there is a pair of  disjoint
  $3$-GDDs of type $g^t u^1$  by the
 Weighting Construction.

 Next we consider even $u$ with $u<2kx= 6x$ for $t=4$. Start from a 4-GDD of type $3^5$ with groups $G_i$, $1\leq i\leq
 5$, where $G_5=\{y_1,y_2,y_3\}$.
 Weight $2x$ to each point of $G_i$ with $1\leq i\leq 4$ and weight even
 weight $w_j,0\leq w_j\leq 4x$, to each point $y_j$ of $G_5$  such that $u=\sum_{j=1}^3 w_j$.  Utilize a pair of disjoint 3-GDDs
  of type $(2x)^4$ or $(2x)^3 w^1$ for even $0\leq w\leq 4x$ and then obtain a pair of  disjoint
  $3$-GDDs of type $g^t u^1$ similarly.

 Finally let  $u$ be even with $u<2kx=2x$ for $t=6$. Start from a $\{4,7\}$-GDD of type
 $3^6 6^1$ with a long group $Y=\{y_1,y_2,\ldots,y_6\}$, in which  precisely one point $y_1$ in $Y$ belongs to blocks of
size $7$. Assign $y_i$ with  $1\leq i\leq 5$ weight 0,  $y_6$ weight
$u$, and each point of the group of size 3 weight $2x$.  Utilize
disjoint pairs of 3-GDDs of types  $(2x)^s$ ($s=3,4,6$) and $(2x)^3
u^1$
 and then obtain a pair of disjoint
  $3$-GDDs of type $(6x)^t u^1$.
This completes the proof. \qed

We summarize the above results on $g\equiv 0$ (mod 6) in a
corollary.
\begin{Corollary}\label{g6x}
The Main Theorem holds for any admissible triple $(g,t,u)$ with
$g\equiv 0$ {\rm (mod $6$)}.
\end{Corollary}

Then the solutions for  $g=2,3, 4$ are ready-made.
\begin{Lemma}\label{g3}%$({\rm \cite{}})$
The Main Theorem holds for any admissible triple $(3,t,u)$.
\end{Lemma}
\proof Since $(3,t,u)$ is admissible, $t$ is even with $t\geq 4$,
$u$ is odd with $u\neq 3$,  and $1\leq u\leq 3(t-1)$.  If $u\geq 5$
and $t\geq 6$, then by Corollary \ref{g6x} there is a pair of
disjoint 3-GDDs of type $6^{t/2} (u-3)^1$.  Apply Corollary
\ref{2gtog} to yield a pair of disjoint $3$-GDDs of type $3^t u^1$.

If $t=4$, then $u=1,5,7,9$. A pair of disjoint 3-GDDs of type $3^4
9^1$ exists by Corollary \ref{u=g(t-1)}. The solutions for $u=1,5,7$
are listed in the appendix.

For $u=1$ and $t=6,8$, let $X=I_3\times I_t$ and ${\cal
G}=\{I_3\times\{i\}:i\in I_t\}\cup\{\infty\}$. First construct on
each $\{j\}\times I_t$ $(j\in I_3)$ a pair of disjoint 3-GDDs of
type $1^{t+1}$. Then form a pair of disjoint ITD$(3,t)$s and delete
their idempotent parallel class. Thus  a pair of disjoint 3-GDDs of
type $3^{t}1^1$ is obtained.

For  $u=1$ and even $t$ with $t\geq 10$, there are pairs of disjoint
3-GDDs of types $3^{t-4} 13^1$ and $3^4 1^1$ by the above arguments.
Consequently a pair of disjoint 3-GDDs of types of $3^{t} 1^1$ is
produced by Filling Construction II.
 \qed

\begin{Lemma}\label{g4}%$({\rm \cite{}})$
The Main Theorem holds for any admissible triple $(4,t,u)$.
\end{Lemma}
\proof Note that $(4,t,u)$ is an admissible triple requires that
$2\leq u\leq 4(t-1)$, $u\neq 4$, $t\equiv 0$ {\rm(mod\ $3$)} and
$u\equiv 0$ {\rm(mod\ $2$)}, or $t\equiv 1$ {\rm(mod\ $3$)} and
$u\equiv 0$ {\rm(mod\ $6$)}, or $t\equiv 2$ {\rm(mod\ $3$)} and
$u\equiv 4$ {\rm(mod\ $6$)}.

Firstly,  when $t\equiv 1$ {\rm(mod\ $3$)} and $u\equiv 0$ {\rm(mod\
$6$)}, or $t\equiv 2$ {\rm(mod\ $3$)} and $u\equiv 4$ {\rm(mod\
$6$)}, or $t\equiv 0$ {\rm(mod\ $3$)} and $u\equiv 2$ {\rm(mod\
$6$)}, let $D=\{1,2,\ldots,2t-1\}\setminus\{t\}$. By Lemma
\ref{dif}, it suffices to show that $D$ can be partitioned into a
set $D_1$ of $(4t-4-u)/6$ difference triples and a set $D_2$
containing a good difference in $Z_{4t}$. This has been done in
Section 4 of \cite{rosa}.

Secondly, let $t\equiv 0$ {\rm(mod\ $3$)}, $u\equiv 0,4$ {\rm(mod\
$6$)}, $ u\geq 6$, and $t\geq 9$. By Corollary \ref{g6x} there is a
pair of disjoint 3-GDDs of type $12^{t/3} (u-4)^1$. A pair of
disjoint 3-GDDs of type $4^4$ also exists by Lemma \ref{intersect}.
Apply Filling Construction I to produce a pair of disjoint $3$-GDDs
of type $4^{t} u^1$.

Finally, we only need to handle $t=3,6$,  $u\equiv 0,4$ {\rm(mod\
$6$)}, and $ u\geq 6$. The case $t=3$ is solved by Corollary
\ref{t3}. There is a pair of disjoint $3$-GDDs of type
$8^{3}(u-4)^1$, so by Corollary \ref{2gtog}, there exists a pair of
disjoint 3-GDD of type $4^6u^1$. \qed

\begin{Lemma}\label{g2}%$({\rm \cite{}})$
The Main Theorem holds for any admissible triple $(2,t,u)$.
\end{Lemma}
\proof By Lemma \ref{g21}, we only need to deal with the admissible
triples $(2,t,u)$ with $t\equiv 0$ (mod 3)  and even $u$ with
$4\leq u\leq 2(t-1)$. If $t\equiv 3$ (mod 6), a pair of disjoint
3-GDDs of type $2^{t}u^1$ is obtained by Corollary
 \ref{t3}. Otherwise, $t\equiv 0$ (mod 6). There exists by Lemma \ref{g4} a pair of
disjoint 3-GDDs of type $4^{t/2}(u-2)^1$. Then the conclusion
follows by Corollary \ref{2gtog}. \qed

To conclude this section we prove that the necessary conditions of
the
 existence of two disjoint 3-GDDs of type $g^t u^1$ for $g\equiv 3$
 (mod $6$) are also sufficient.

\begin{Lemma}\label{g6k+3}%$({\rm \cite{}})$
The Main Theorem holds for any admissible triple $(g,t,u)$ with
$g\equiv 3$ {\rm (mod\ $6$)}.
 \end{Lemma}
\proof Since $g\equiv 3$ {\rm (mod\ $6$)} and $(g,t,u)$ is
admissible, $t$ must be even with $t\geq 4$, $u$ be odd,  and $u\leq
g(t-1).$ Let $(X,{\cal A})$ be a KTS$(g)$, where ${\cal A}$ can be
resolved into $(g-1)/2$ parallel classes
$P_1,P_2,\ldots,P_{(g-1)/2}$. Choose integers $u_i$, $1\leq i\leq
(g-1)/2,$ such that $u_1$ is odd, $1\leq u_1\leq 3(t-1)$ and for
each $2\leq i\leq (g-1)/2,$ $u_i$ is even, $0\leq u_i\leq 2(t-1).$
Let $U_1, U_2,\ldots, U_{(g-1)/2}$ be pairwise disjoint sets with
$|U_i|=u_i$ and let $U=\cup_{i=1}^{(g-1)/2}U_i$. The desired two
disjoint 3-GDDs will be constructed on the set $Y=(X\times I_t)\cup
U$ with  group set ${\cal G}=\{X\times\{i\}:i\in I_t\}\cup\{U\}$.

For each block $B=\{x,y,z\}\in P_1$, there is a pair of disjoint
3-GDDs $(X_B,{\cal G}_B, {\cal A}_B^1)$ and $(X_B,{\cal G}_B, {\cal
A}_B^2)$ of type $3^t {u_1}^1$ by Lemmas \ref{intersect} and
\ref{g3}, where $X_B=(B\times I_t)\cup U_1$ and ${\cal
G}_B=\{B\times\{i\}:i\in I_t\}\cup\{U_1\}$.

For each block $B=\{x,y,z\}\in P_i$, $2\leq i\leq (g-1)/2,$ there is
a pair of 3-GDDs of type $t^3 {u_i}^1$ with no block in common but a
common parallel class $P=\{B\times \{i\}:i\in I_t\}$ of $B\times
I_t$ by Corollary \ref{t3i}. Deleting the common parallel class $P$
yields
 two disjoint block sets ${\cal A}_B^1$ and ${\cal A}_B^2$.

 For $i=1,2$, let ${\cal B}_i=\cup_{B\in P_j,1\leq j\leq(g-1)/2}{\cal
 A}_B^i.$ It can be checked that $(Y,{\cal G},{\cal B}_1)$ and $(Y,{\cal G},{\cal
 B}_2)$ form a  pair of
disjoint $3$-GDDs of type $g^{t} u^1$. \qed

%\section{The case $g\equiv 1,2$ $({\rm mod}$ $3)$}

\section{Further constructions}

In this section, we shall go a step further to employ cyclic partial
S$(2,3,v)$s to construct a pair of disjoint 3-GDDs.

\begin{Lemma}\label{u>t,geven}
Suppose that $g$ is an even integer and there is  a cyclic partial
S$(2,3,g)$ which contains a starter block having a good difference
and whose leave is $r$-regular. Let $t\geq 4$ and $t\neq 6,10$,
$0\leq m\leq t-1$, and $0\leq v\leq 2(t-1)$ such that a pair of
disjoint $3$-GDDs of type $2^tv^1$ exists. Then there is a pair of
disjoint $3$-GDDs of type $g^t ((r-1)(t-1)+6m+v)^1.$
\end{Lemma}
\proof Let $G=\{\infty_1,\infty_2,\ldots,\infty_v\}$, $X=(Z_g\times
I_t)\cup G$, and ${\cal G}=\{Z_g\times\{i\} :i\in I_t\}\cup\{G\}$.
 For $D\subseteq Z_g$, $x\in Z_g$, denote
$D+x=\{d+x:d\in D\}$ and $dev(D)=\{D+x:x\in Z_g\}$. For
$\Omega\subseteq Z_g\times I_t$, $x\in Z_g$, denote
$\Omega+x=\{(d+x,i):(d,i)\in \Omega\}$ and
$dev(\Omega)=\{\Omega+x:x\in Z_g\}$.

Let $S_1,S_2,\ldots,S_n$ be the starter blocks of a cyclic partial
S$(2,3,g)$ on $Z_g$, whose $r$-regular leave is $L$. Further suppose
that $S_1$ contains a good difference. Clearly, $g/2$ appears as a
difference in $L$  but not in $S_1$.  Let
$L_1=\bigcup_{\{a,b\}\subseteq S_1} dev(\{a,b\})$. By Lemma
\ref{goodd} and noting that $S_1$ contains a good difference, $L$
has a 1-factorization with 1-factors $F_1,F_2,\ldots,F_r$ and $ L_1$
has also a 1-factorization with $H_1,H_2,\ldots,H_6$, as 1-factors.

First for each pair $P\in F_1$, we can construct by the assumption
on $(P\times I_t)\cup G$ a pair of disjoint 3-GDDs of type $2^t v^1$
with group set $\{P\times \{i\}:i\in I_t\}\cup\{G\}$ and two
disjoint block sets ${\cal C}_{P}^0$ and ${\cal C}_{P}^1$. Set
${\cal C}^s=\bigcup_{P\in F_1}{\cal C}_{P}^s$ for $s=0,1$. (The
other 1-factors are left for later use.)

Next we employ the starter block $S_1$. By Lemma \ref{itd}, for
$t\geq 4$ and $t\neq 6, 10,$ there is a pair of disjoint
RITD$(3,t)$s on $S_1\times I_t$  with group set $\{\{x\}\times
I_t:x\in S_1\}$. Let $P_0^s,P_1^s,\ldots,P_{t-1}^s$ ($s=0,1$) be
their parallel classes, where $P_0^s$ be the idempotent one. By
deleting $m+1$ parallel classes,  $P_k^s$, $0\leq k\leq m$, we
obtain two disjoint partial 3-GDDs with block sets ${\cal B}_1^0$
and ${\cal B}_1^1$.

Then we employ the starter block $S_i$ $(i\neq 1)$. For each $2\leq
i\leq n$, construct on $S_i\times I_t$ two disjoint ITD$(3,t)$s with
group set $\{\{x\}\times I_t:x\in S_i\}$. Delete the idempotent
parallel class to form two disjoint block sets ${\cal B}_i^0$ and
${\cal B}_i^1$.

After that, for  $s=0,1$, define ${\cal B}^s=\bigcup_{1\leq i\leq
n}dev({\cal B}_i^s)$ and ${\cal A}^s={\cal B}^s\cup{\cal C}^s.$ One
can check that $(X,{\cal G},{\cal A}^0)$ and $(X,{\cal G},{\cal
A}^1)$ form two disjoint partial 3-GDDs of type $g^t v^1$ with
leaves ${\cal L}^0$ and ${\cal L}^1$. If $(r-1)(t-1)+6m=0$, then
${\cal L}^s$ is empty and we do have obtained a pair of disjoint
$3$-GDDs of type $g^t ((r-1)(t-1)+6m+v)^1.$ So we assume that $r\geq
2$ or $m\geq 1$. By the previous construction, for $s=0,1$, ${\cal
L}^s$
 consists of two parts ${\cal L}_1^s$ and ${\cal L}_2^s$, where  ${\cal
L}_1^0={\cal L}_1^1=\{\{(a,i),(b,j)\}:\{a,b\}\in L\setminus
F_1,i\neq j\in I_t\}$, and ${\cal L}_2^s$ contains all the pairs in
$\bigcup_{k=1}^{m}dev(P_k^s).$

Finally we partition each ${\cal L}^s$ into $(r-1)(t-1)+6m$ disjoint
1-factors of $Z_g\times I_t$ to complete the proof. For $\{a,b\}\in
L\setminus F_1$ and $1\leq i\leq t-1$, take
$f_{ab}^i=\{\{(a,j),(b,j+i)\}:0\leq j\leq t-1\}$. Then we have $t-1$
disjoint 1-factors of $\{a,b\}\times I_t.$  For $\{a,b\}\in L_1$ and
$Q= dev(P_{k}^s)$ ($1\leq k\leq m$ and $s=0,1$), take
$f_{ab}^Q=\{\{(a,l),(b,u)\}:\{(a,l),(b,u),(c,w)\}\in Q\}$. Thus we
have $m$ disjoint 1-factors of $\{a,b\}\times I_t$ for each $s=0,1$,
which for convenience we also denote in sequence by
$f_{ab}^{s1},f_{ab}^{s2},\ldots,f_{ab}^{sm}$. Define
$$D_{ij}=\bigcup_{\{a,b\}\in
F_j}\{\{\alpha,\beta\}:\{\alpha,\beta\}\in
 f_{ab}^i\}, {\rm where} \ 1\leq i\leq t-1\ {\rm and}\ 2\leq j\leq
 r,$$
 $$E_{kl}^s=\bigcup_{\{a,b\}\in H_l}\{\{\alpha,\beta\}:\{\alpha,\beta\}\in
 f_{ab}^{sk}\}, {\rm where}\ 1\leq k\leq m\ {\rm and}\ 1\leq l\leq
6.$$ It is readily checked that the union of these $D_{ij}$'s and
$E_{kl}^s$'s equals ${\cal L}^s$, forming $(r-1)(t-1)+6m$ disjoint
1-factors of $Z_g\times I_t$. Obviously the number of these
1-factors is greater than 2 when $t\geq 4$ and $r\geq 2 $ or $m\geq
1$, so we can arrange them such that Lemma \ref{1-factor}  can be
applied to form a pair of disjoint $3$-GDDs of type $g^t
((r-1)(t-1)+6m+v)^1.$
 \qed

 For any integer $g\geq 2$, there is a trivial cyclic S$(2,3,g)$
 (with no starter block) whose leave is $(g-1)$-regular. Then in a
 similar but simpler  procedure than the proof of Lemma \ref{u>t,geven},
 we have an analogous result (the details of the proof are omitted).
 \begin{Lemma}\label{r=g-1,geven}
Suppose that $g$ is an even integer. Let $t\geq 4$, $t\neq 6,10$,
$0\leq m\leq t-1$, and $0\leq v\leq 2(t-1)$ such that a pair of
disjoint $3$-GDDs of type $2^tv^1$ exists. Then there is a pair of
disjoint $3$-GDDs of type $g^t ((g-2)(t-1)+v)^1.$
\end{Lemma}

\begin{Lemma}\label{connected}{\rm (\cite{stong})}
Suppose that $\Gamma$ is an abelian group of even order and
$S\subseteq \Gamma\setminus\{0\}$. Let $G(\Gamma,S)$ be the graph
with vertex set $\Gamma$ and whose edge set is
$\{\{x,x+s\}:x\in\Gamma,s\in S\}$. Then $G(\Gamma,S)$ has a
$1$-factorization whenever it is connected.
\end{Lemma}

\begin{Lemma}\label{u>t,godd}
Suppose that there is  a cyclic partial S$(2,3,g)$ whose leave is
$r$-regular with $r<g-1$. Let $t\geq 4$ be even, $0\leq m\leq t-1$,
and $1\leq v\leq t-1$ such that a pair of disjoint $3$-GDDs of type
$1^tv^1$ exists. Then there is a pair of disjoint $3$-GDDs of type
$g^t (r(t-1)+6m+v)^1.$
\end{Lemma}
\proof Let $G=\{\infty_1,\infty_2,\ldots,\infty_v\},$ $X=(Z_g\times
I_t)\cup G$, and ${\cal G}=\{Z_g\times \{i\}:i\in I_t\}\cup\{G\}$.
We first construct two disjoint partial $3$-GDDs of type $g^t v^1$
on $X$ with group set ${\cal G}$. Then we partition their leaves
into $r(t-1)+6m$ disjoint 1-factors.  For $D\subseteq Z_g$,
$\Omega\subseteq Z_g\times I_t$, and $x\in Z_g$, we use the
notations $D+x$, $\Omega+x$, $dev(D)$, and $dev(\Omega)$ as in Lemma
\ref{u>t,geven}.

By the assumption, for each $i\in Z_g$, there is a  pair of 3-GDDs
of type $1^t v^1$  on $(\{i\}\times I_t)\cup G$ with $G$ as the long
group and disjoint block sets ${\cal D}_i^0$ and ${\cal D}_i^1$. For
$s=0,1,$ set ${\cal D}^s=\cup_{i\in Z_g}{\cal D}_i^s$.

Let $S_1,S_2,\ldots,S_n$ be the starter blocks of the cyclic partial
S$(2,3,g)$ on $Z_g$, whose $r$-regular leave is $L$. For each $2\leq
i\leq n$,  construct on $S_i\times I_t$ two disjoint ITD$(3,t)$s
with group set $\{\{x\}\times I_t:x\in S_i\}$ and delete the
idempotent parallel class to form two disjoint block sets ${\cal
C}_i^0$ and ${\cal C}_i^1$.

 Next we handle $S_1$. Let $S_1=\{a,b,c\}$. If $m=0$, we deal with $S_1$ as
$S_i$. So suppose $m\geq 1$. For $t\geq 6$ and $t\neq 12,$ there is
an RITD$(3,t/2)$ on $S_1\times\{2k:0\leq k\leq t/2-1\}$  with group
set $\{\{x\}\times\{2k:0\leq k\leq t/2-1\}:x\in S_1\}$ and $t/2$
parallel classes $P_1,P_2,\ldots,P_{t/2}$, where
$P_1=\{S_1\times\{2k\}:0\leq k\leq t/2-1\}$. Define $M= (t-m+1)/2$
if $m$ is odd,  or $M= (t-m+2)/2$ if $ m$ is even. We proceed with
$M$ parallel classes as follows:

 Take any block $B=\{(a,2i),(b,2j),(c,2k)\}\in P_l$, $l=1$ if $m$ is odd, or $l=1,2$
if $m$ is even. For  $s=0,1$, form a partial 3-GDD of type $2^3$
with group set $\{\{a\}\times
\{2i+2s,2i+2s+1\},\{b\}\times\{2j,2j+1\},\{c\}\times\{2k,2k+1\}\}$
and block set ${\cal A}_B^s$, where
\begin{equation}\label{*}
{\cal
A}_B^{s}=\{\{(a,2i+2s),(b,2j),(c,2k)\},\{(a,2i+2s+1),(b,2j+1),(c,2k+1)\}\},
\end{equation}
  and the second components  are modulo $t$. %Denote by
%$L_B^s$ the leave of ${\cal A}_B^s$.

For any block $B=\{(a,2i),(b,2j),(c,2k)\in P_l$, $2\leq l\leq M$ if
$m$ is odd, or $3\leq l\leq M$ if $m$ is even, take a 3-GDD with
group set $\{\{a\}\times
\{2i+2s,2i+2s+1\},\{b\}\times\{2j,2j+1\},\{c\}\times\{2k,2k+1\}\}$
and block set ${\cal A}_B^s$, where  $s=0,1$.

 For $s=0,1,$ define ${\cal
 C}_1^s=\bigcup_{B\in P_l,1\leq l\leq M}\{dev(A):A\in {\cal A}_B^s\}$.
 Then by defining ${\cal C}^s=\bigcup_{i=1}^{n}{\cal C}_i^s$ and
 ${\cal B}^s={\cal D}^s\bigcup{\cal C}^s$, we produce two disjoint partial
3-GDDs of type $g^t v^1$ $(X,{\cal G},{\cal B}^0)$ and $(X,{\cal
G},{\cal B}^1)$. Denote their leaves by ${\cal L}_0$ and ${\cal
L}_1$, respectively. By the construction, ${\cal L}_s$ $(s=0,1)$
consists of at most three parts. We partition the pairs in the leave
into $r(t-1)+6m$ disjoint 1-factors of $Z_g\times I_t$ to complete
the proof for $t\geq 6$ and $t\neq 12$.

Part I: For $s=0,1,$ $l=1$ if $m$ is odd, or $l=1,2$ if $m$ is even,
observe that we take a partial 3-GDD as in the expression (\ref{*})
for each block $B=\{(a,2i),(b,2j),(c,2k)\}$  of $P_l$, leading to
the leave ${\cal L}_l^{s}={\cal L}_{l0}^{s}\cup {\cal
L}_{l1}^{s}\cup {\cal L}_{l2}^{s}$ with
\begin{eqnarray*}
&&{\cal L}_{l0}^{s}=\bigcup_{B\in
P_l}(dev(\{(a,2i+2s),(b,2j+1)\})\cup
dev(\{(a,2i+2s+1),(b,2j)\})),\\
&&{\cal L}_{l1}^{s}=\bigcup_{B\in
P_l}(dev(\{(a,2i+2s),(c,2k+1)\})\cup
dev(\{(a,2i+2s+1),(c,2k)\})),\\
&&{\cal L}_{l2}^{s}=\bigcup_{B\in P_l}(dev(\{(b,2j),(c,2k+1)\})\cup
dev(\{(b,2j+1),(c,2k)\})).
\end{eqnarray*}
Observe that the second components of each pair in  ${\cal
L}_{li}^{s}$ ($i=0,1,2$) are not equivalent modulo 2.  So the graph
${\cal L}_{li}^{s}$  consists of some cycles of even length. Thus
each cycle has a 1-factorization with two 1-factors. By collecting
the 1-factors corresponding to all the connected cycles of ${\cal
L}_{li}^{s}$, we obtain  two 1-factors of $Z_g\times I_t$, say
$F_{l,2i}^{s}$ and $F_{l,2i+1}^{s}$. Furthermore, $F_{l,p}^0\cap
F_{l,p+2}^1=\emptyset,$ where $p\in I_6$ and $p+2$ is reduced to
$I_6$.  Now for fixed $s$ we have six 1-factors of $Z_g\times I_t$
for odd $m$ or twelve 1-factors for even $m$.

Part II: This part of leave exists only if $m\geq 3$. For $s=0,1,$
and $M+1\leq l\leq t/2$, observe that we do not use any block in
$P_l$, which leads to leave ${\cal L}_l^s$ described below. For each
$B=\{(a,2i),(b,2j),(c,2k)\}\in P_l$, ${\cal L}_l^s$ contains the
pairs in the 2-GDD with group set $\{\{a\}\times
\{2i+2s,2i+2s+1\},\{b\}\times\{2j,2j+1\},\{c\}\times\{2k,2k+1\}\}$.
By similar arguments, ${\cal L}_l^s$ can be partitioned into twelve
disjoint 1-factors  of $Z_g\times I_t$ and we obtain $K$ 1-factors
altogether, say, $G_0^s,G_1^s,\ldots, G_{K-1}^s$, where $K=6(m-1)$
for odd $m$ or $K=6(m-2)$  for even $m$. Furthermore, we can arrange
them such that $G_i^0\cap G_i^1=\emptyset$ holds for all $0\leq
i\leq K-1$.

Part III: This part of leave exists only if $r\neq 0$. We consider
the leave $L$ of the cyclic partial S$(2,3,g)$. Observe that
$dev(P)$ is a 2-regular graph consisting of some cycles for any pair
$P\in L$. For each connected component $C$, the set
$\{\{(u,i),(w,j)\}:\{u,w\}\in C,i\neq j\in I_t\}$  can be
1-factorized by Lemma \ref{connected} (taking $\Gamma=\{(i\ {\rm
mod}\ |C|,i\ {\rm mod}\ t):0\leq i\leq {\rm lcm} (|C|,t)\}$ and
$S=\{0\}\times (Z_t\setminus\{0\})$). Thus $r(t-1)$ 1-factors of
$Z_g\times I_t$ are obtained when taking $P$ all over the
$r$-regular leave $L$. These 1-factors,
$H_0,H_1,\ldots,H_{r(t-1)-1}$, are all contained in both ${\cal
L}_0$ and ${\cal L}_1$ and certainly $H_i\cap H_{i+1}=\emptyset$.

So we obtain $r(t-1)+6m$ disjoint 1-factors altogether. By Lemma
\ref{1-factor}, there is a pair of disjoint $3$-GDDs of type $g^t
(r(t-1)+6m+v)^1$ for $t\geq 6$ and $t\neq 12$.

If $t=4$, we can utilize on $S_1\times I_4$ an RITD$(3,4)$ with the
idempotent parallel class omitted and further empty some parallel
classes. If $t=12$,  we use on $S_1\times \{3k:0\leq k\leq 3\}$ an
RITD$(3,4)$ with the idempotent parallel class omitted. And then
deal with its four parallel classes by two ways. Choose appropriate
number of parallel classes to construct for each $s=0,1$ an RTD(3,3)
with groups $\{a\}\times\{3i+3s,3i+3s+1,3i+3s+2\}$,
$\{b\}\times\{3j+3s,3j+3s+1,3j+3s+2\}$, and
$\{c\}\times\{3k+3s,3k+3s+1,3k+3s+2\}$, where
$\{(a,3i),(b,3j),(c,3k)\}$ is any block of the chosen parallel
classes. And for each block of the remaining parallel classes of the
RITD$(3,4)$, also take RTD(3,3) similarly but delete some parallel
classes of this RTD. Then in a very similar way, a pair of disjoint
$3$-GDDs of type $g^t (r(t-1)+6m+v)^1$ is constructed. This
completes the proof. \qed

Parallel to Lemma \ref{r=g-1,geven}, the following result also
holds.
 \begin{Lemma}\label{r=g-1,godd}
Suppose that $g$ is a positive integer. Let $t\geq 4$ be even,
$0\leq m\leq t-1$, and $1\leq v\leq t-1$ such that a pair of
disjoint $3$-GDDs of type $1^tv^1$ exists. Then there is a pair of
disjoint $3$-GDDs of type $g^t ((g-1)(t-1)+v)^1.$
\end{Lemma}

{\Lemma\label{most}  Let $(g,t,u)$ be any admissible triple with
$g>5$ and $t\geq 4$. Then there exists a pair of disjoint $3$-GDDs
of type $g^t u^1$ whenever one of the following conditions meets:
\begin{enumerate}
\item[{\rm (1)}]  $g\equiv 2,8$ {\rm(mod\ $24$)} if $t\neq 6,10;$

\item[{\rm (2)}] $g\equiv 14,20$ {\rm(mod\ $24$)} and  $u\geq 6(t-1)$ if $t\neq 6,10;$

\item[{\rm (3)}]  $g\equiv 4$ {\rm(mod\ $6$)} and $u\geq 2(t-1)$  if $t\neq 6,10;$

\item[{\rm (4)}] $g\equiv 1$ {\rm(mod\ $6$);}

\item[{\rm (5)}] $g\equiv 5$ {\rm(mod\ $6$)} and $u> 4(t-1);$

\item[{\rm(6)}] If $t=6,10$, then $u>t-1$ for $g\equiv 2,8$ {\rm(mod\ $24$),}
or $u>7(t-1)$ for $g\equiv 14,20$ {\rm(mod\ $24$),} or $u>3(t-1)$
for $g\equiv 4$ {\rm(mod\ $6$).}
\end{enumerate}
}

\proof Suppose that $g=6k+s$, where $k\geq 1$ and $1\leq s\leq 6.$
Let $r'=7$ if $s=2$ and $k\equiv 2,3$ {\rm(mod\ $4$)}, or $r'=s-1$
otherwise.

  For any admissible $(g,t,u)$ with $g\equiv 2,4$ {\rm(mod\ $6$)},
$t\geq 4$, $t\neq 6,10$, and $u\geq
  (r'-1)(t-1)$,  first take $0\leq x<6,$ $x\equiv u-(r'-1)(t-1)$ (mod 6) ($x$ must be even) and next choose  $ r\equiv r'$ {\rm(mod\
 $6$)} and $0\leq u-(r-1)(t-1)-x=6m\leq 6(t-1),$ then $u=(r-1)(t-1)+6m+x$ and $r\leq g-1$. By Lemma \ref{cps}, there is a cyclic partial S$(2,3,g)$
 with an $r$-regular leave. Moreover, if  $r<g-1$, there is a starter block
containing a good difference. And we can check that $(2,t,x)$ is
 an admissible triple and then obtain a pair of disjoint 3-GDDs of type
 $2^t x^1$ by Lemma \ref{g2}. Consequently there is a pair of disjoint 3-GDDs of type
 $g^t u^1$ by Lemma \ref{u>t,geven}. If $r=g-1$, then $(2,t,6m+x)$ is admissible and  a pair of disjoint 3-GDDs of type
 $2^t (6m+x)^1$ also exists. So the conclusion follows by Lemma \ref{r=g-1,geven}. This
 handles (1)-(3).

  For any admissible $(g,t,u)$ with $g\equiv 1,5$ {\rm(mod\ $6$)} (or $g\equiv 2,4$ {\rm(mod\ $6$)} and $t=6,10$) and  $u> r'(t-1)$,
  first take $0\leq x<6,$ $x\equiv u-r'(t-1)$ (mod 6) ($x$ must be odd) and next choose  $ r\equiv r'$ {\rm(mod\
 $6$)} and $0\leq u-r(t-1)-x=6m\leq 6(t-1),$ then $u=r(t-1)+6m+x$ and $r\leq g-1$.
 By Lemma \ref{cps}, there is a cyclic partial S$(2,3,g)$
 with an $r$-regular leave. It can be checked that
 $(1,t,x)$ (if $r<g-1$) or $(1,t, 6m+x)$ (if $r=g-1$) is
 an admissible triple, so there is a pair of disjoint 3-GDDs of type
 $1^t x^1$ or $1^t (6m+x)^1$ by Lemma \ref{g1}.  Consequently there is a pair of disjoint 3-GDDs of type
 $g^t u^1$ by Lemma \ref{u>t,godd} or \ref{r=g-1,godd}. This proves (4)-(6).
\qed

%Remark: By a similar method as in the above lemma, we can also solve
%some cases of $g\equiv 3$ {\rm(mod\ $6$)} when $u>2(t-1)$, such as
%$u\equiv 1,3,5$ (mod 6) if $t\equiv 0$ (mod 6), $u\equiv 1$ (mod 6)
%if $t\equiv 2$ (mod 6), and $u\equiv 3$ (mod 6) if $t\equiv 4$ (mod
%6).

\section{The case $g\equiv 2,4$ $({\rm mod}$
$6)$}

We handle the remaining cases when $g\equiv 2,4$ $({\rm mod}$ $6)$
in this section.

{\Lemma \label{g6k+4}  The Main Theorem holds for any admissible
triple  $(g,t,u)$  with $g\equiv 4$ {\rm(mod $6$)}. }

\proof  By Lemma \ref{most}, we need only to consider admissible
triples with $u<2(t-1)$ if $t\neq 6,10$ and $u\leq 3(t-1)$ if
$t=6,10$. Let $g=6n+4$. The case $n=0$ or $t=3$ is solved by Lemma
\ref{g4} and Corollary \ref{t3} respectively. So suppose that $n\geq
1$ and $t\geq 4$. Since $(g,t,u)$ is admissible, either $u\equiv 0$
(mod 2) if $t\equiv 0$ (mod 3), or  $u\equiv 0$ (mod 6) if $t\equiv
1$ (mod 3), or  $u\equiv 4$ (mod 6) if $t\equiv 2$ (mod 3).  We
distinguish all the possible cases.

Case 1: $n\geq 3$ and $u\leq 3(t-1)$.  There is a 3-GDD of type $6^n
4^1$ by Lemma \ref{2.3}. There are pairs of disjoint 3-GDDs of types
$6^t u^1$ and $4^tu^1$ by Corollary \ref{g6x} and Lemma \ref{g4}. So
a pair of disjoint 3-GDDs of type $(6n+4)^t u^1$ is obtained by
Construction \ref{glarge}.

Case 2: $n=2$ and $u\leq 3(t-1)$. There  is a 3-GDD of type
 $4^4$ by Lemma \ref{2.3}.  There is a  pair of disjoint 3-GDDs of
type $4^t u^1$ by Lemma \ref{g4}. So there exists a pair of disjoint
3-GDDs of type $16^t u^1$ by Construction \ref{glarge}.

Case 3: $n=1,$ $t\equiv 2$ (mod 3), and $u< 2(t-1)$. Then $g=10$ and
$u\equiv 4$ (mod 6). First Lemma \ref{u>g} solves  such cases with
$u\geq 2g+2=22$, leaving $u=4$ if $t\leq 8$ or $u=4,16$ if $t\geq11$
to be settled. Next utilize Lemma \ref{dif} to deal with $t=5$ and
$u=4$ by taking on $Z_{50}$ the difference triples $\{1,23,24\}$,
$\{4,18,22\}$, $\{6,7,13\}$, $\{8,11,19\}$, $\{9,12,21\}$ and
$\{2,14,16\}$. Finally  for $t=8$ and $u=4$,  or $t\geq 11$ and
$u=4,16$, the Filling Construction II works by filling a pair of
disjoint 3-GDDs of type $10^{t-3}(30+u)^1$  with such pair of type
$10^3 u^1$.

Case 4: $n=1,$ $t\equiv 0,1$ (mod 3), and $u< 2(t-1)$.  There  is a
3-GDD of type $2^{3}4^1$
  and disjoint pairs of 3-GDDs of types $2^tu^1 $ and $4^tu^1$ exist by Lemmas \ref{g4} and \ref{g2}. So we produce
a pair of disjoint 3-GDDs of type $10^tu^1$ by Construction
\ref{glarge}.

Case 5: $n=1$, $t=6,10$, and $2(t-1)\leq u\leq 3(t-1)$. Then $u\geq
10$ if $t=6$. So there exists a pair of disjoint 3-GDDs of type
$10^6 u^1$ by Corollary \ref{2gtog} since there is a pair of
disjoint 3-GDDs of type $20^{3} (u-10)^1$ by Corollary \ref{t3}. If
 $t=10$, then $u=18,24$. Thus a pair of disjoint 3-GDDs of type $10^{10}
u^1$ exists by Lemma \ref{u>g}.
 \qed

{\Lemma \label{g1420}  The Main Theorem holds for any admissible
triple $(g,t,u)$  with $g=14,20$. }

\proof For $g=14,20$, the case $t\equiv 3$ (mod 6) has been solved
by Corollary \ref{t3}, so let
 $t\not\equiv 3$ (mod 6). If $t\geq 6$ is even and $u> g$, a pair of disjoint
3-GDDs of type $g^tu^1$ can be obtained by Corollary \ref{2gtog}
since a pair of disjoint 3-GDDs of type $(2g)^{t/2}(u-g)^1$ exists
by Lemma \ref{g6k+4}. Thus by Lemma \ref{most} we need only to
consider $u<g$ if $t\geq 6$ is even and $u<6(t-1)$ if $t=4$ or
$t\equiv 1,5$ (mod 6).  Since $(g,t,u)$ is admissible, either
$u\equiv 0$ (mod 2) if $t\equiv 0$ (mod 3), or $u\equiv 0$ (mod 6)
if $t\equiv 1$ (mod 3), or  $u\equiv 2$ (mod 6) if $t\equiv 2$ (mod
3).

(1) $g=14$.

Case 1:  $t\geq 5$ and $u<14$. Then $u\leq2(t-1)$ (noting that
$(g,t,u)$ is admissible) and there exist a 3-GDD of type $2^7$ and a
pair of disjoint 3-GDDs of type $2^tu^1$ by Lemma \ref{g2}, yielding
a pair of disjoint 3-GDDs of type $14^tu^1$ by Construction
\ref{glarge}.

Case 2: $t=4,6,7$ and $u<6(t-1)$, or $t=5$ and $14\leq u<6(t-1)=24$.
Employ the Weighting Construction. Start from a TD$(t+1,7)$. Assign
weight 2 to each point of the first $t$ groups and then assign
appropriate weight $w$ to the point of the last group, where
$w\equiv 0$ (mod 2) if $t=6$, or  $w\equiv 0$ (mod 6) if
$t\in\{4,7\}$, or  $w\equiv 2$ (mod 6) if $t=5.$

Case 3: $t\equiv 1,5$ (mod 6), $t\geq 9$, and $u<6(t-1)$.   First
Lemma \ref{u>g} solves such cases with $u\geq 2g+2=30$, leaving
$u\leq 28$ to be settled. Then fill a pair of disjoint 3-GDDs of
type $14^3 u^1$ in that of type $14^{t-3}(42+u)^1$ to obtain a pair
of disjoint 3-GDDs of type $14^tu^1$.

(2) $g=20$.

Case 1: $t\equiv 1,5$ (mod 6), $t\geq 11$ and $u<6(t-1)$. Similarly
Lemma \ref{u>g} solves such cases with $u\geq 2g+2=42$. For $u\leq
40$, fill in the long group of a pair of disjoint 3-GDDs of type
$20^{t-3}(60+u)^1$ with that of type $20^3u^1$ to produce the
desired pair of type $20^tu^1$.

Case 2: even $t\geq 10$ and $u<20$, or $t=5$ and $u<6(t-1)=24$. If
$t=5$ and $u=14$, employ Lemma \ref{dif} on $Z_{100}$ by taking
difference triples $\{1,2,3\}$, $\{4,7,11\}$, $\{6,8,14\}$,
$\{9,12,21\}$, $\{13,16,29\}$, $\{17,19,36\}$, $\{18,23,41\}$,
$\{22,24,46\}$, $\{26,27,47\}$, $\{28,33,39\}$, and $\{31,32,37\}$.
If $t\neq 5$ or $u\neq 14$, then $u\leq 2(t-1)$. So these cases can
be solved similarly to the Case 1 of $g=14$, using a 3-GDD of type
$2^{10}$ instead of $2^7$.

Case 3: $t=4$ and $u<6(t-1)=18$, or $t=6$ and $u<20$. Then $u\leq
4(t-1)$ and we can apply Construction \ref{glarge} to a 3-GDD of
type $4^38^1$. A pair of disjoint 3-GDDs of type $4^tu^1$ exist by
Lemmas \ref{g4}. If $t=4,$ or $t=6$ and $u\geq 6$, a pair of
disjoint 3-GDDs of type $8^tu^1$  exists by Lemma \ref{most}. And if
$t=6$ and $u=2,4$, a pair of disjoint 3-GDDs of types $8^tu^1$ also
exists since a  3-GDD of type $2^4$ and a pair of disjoint 3-GDDs of
type $2^tu^1$ exist. Thus
  Construction \ref{glarge} gives a pair of disjoint 3-GDDs of type $20^tu^1$.

Case 4: $t=8$ and $u<20$. Then $u=2,8,14.$ Similar to Case 1, fill
in the long group of a pair of disjoint 3-GDDs of type
$20^{5}(60+u)^1$ with that of type $20^3u^1$ to produce the desired
pair of type $20^8u^1$.

 Case 5: $t=7$ and $u<6(t-1)=36$.
Then $u=6,12,18,24,30.$ As in Case 3, we can handle $u\leq 24$. The
last case $u=30$ is treated as follows.

Let $(X, {\cal G},{\cal B})$ be a $\{2,3\}$-GDD of type $4^5$, which
is obtained by deleting a group of a 3-GDD of type $4^6$. So the
blocks of size 2 of ${\cal B}$ is partitioned into four parallel
classes of $X$. Let $U=\{\infty_1,\infty_2,\ldots,\infty_6\}$,
$Y=(X\times I_7)\cup U$, and ${\cal H}=\{X\times\{i\}:i\in
I_7\}\cup\{U\}$. For each $B\in {\cal B}$ and $|B|=3$, construct on
$B\times I_7$ a pair of disjoint RITD(3,7)s (but deleting the
idempotent parallel class) with group set $\{\{x\}\times I_7:x\in
B\}$ and block sets ${\cal A}_B^1$ and ${\cal A}_B^2$. For each
$G\in{\cal G}$, construct on $(G\times I_7)\cup U$ a pair of
disjoint 3-GDDs of type $4^76^1$ with group set $\{\{x\}\times
I_7:x\in G\}\cup\{U\}$ and block sets ${\cal C}_G^1$ and ${\cal
C}_G^2$. Set ${\cal C}^i=(\cup_{B\in{\cal B},|B|=3}{\cal
A}_B^i)\cup(\cup_{G\in{\cal G}}{\cal C}_G^i)$ where $i=1,2$. Then
$(Y,{\cal H},{\cal C}^1)$ and $(Y,{\cal H},{\cal C}^2)$ form a pair
of disjoint partial 3-GDDs of type $20^76^1$. Their common leave is
$\{((x,i),(y,j)):\{x,y\}\in {\cal B},i,j\in I_7,i\neq j\}$. Noting
that the pairs of ${\cal B}$ is partitioned into four parallel
classes, we can partition the leave into $6\times 4=24$ disjoint
1-factors of $X\times I_7$. Hence there is a pair of disjoint 3-GDDs
of type $20^7 30^1$ by Lemma \ref{1-factor}. \qed

{\Lemma \label{g6k+2}  The Main Theorem holds for any admissible
triple $(g,t,u)$  with  $g\equiv 2$ {\rm(mod $6$).}
 }

 \proof  By Lemmas \ref{most} and \ref{g1420}, for  $g\equiv 2,8$ (mod 24), we need
only to consider $t=6,10$ and $u\leq t-1$. For $g\equiv 14,20$ (mod
24), we need only to consider $g\geq 38$ and $u<6(t-1)$, further
$u\leq 7(t-1)$ if $t=6,10$.  The possible cases are listed as
follows:

Case 1: $g\equiv 2,8$ (mod 24), $t=6,10$,  and $u\leq t-1$. Let
$g=6n+2$. The case $n=0$ is solved by Lemma \ref{g2}. So let $n\geq
1$. Since there are a 3-GDD of type $2^{3n+1}$ and a pair of
disjoint $3$-GDDs of type $2^tu^1$ by Lemmas \ref{2.3} and \ref{g2},
there is a pair of disjoint $3$-GDDs of type $(6n+2)^tu^1$ by
Construction \ref{glarge}.

Case 2: $g\equiv 14,20$ (mod 24), $g\geq 38$,  and $u<6(t-1)$. Let
$g=6l+8$, where $l\geq 5$.  There exists a pair of disjoint 3-GDDs
of type $(6l+8)^t 8^1$ by Construction \ref{glarge} since there are
a 3-GDD of type $6^{l}8^1$ and disjoint pairs of 3-GDDs of types
$6^{t} u^1$ and $8^tu^1$ by Corollary \ref{g6x} and Lemma \ref{most}
or Case 1 of the proof.

Case 3: $g\equiv 14$ (mod 24), $t=6, 10$, and $6(t-1)\leq u\leq
7(t-1)$, where $m\geq 1$. Employ a 3-GDD of type $8^{3m} 14^1$ and
 disjoint pairs of 3-GDDs of types
 $8^t u^1$ and  $14^t u^1$ (whose existence is assured by Case 1 and Lemma \ref{g1420}). Then we obtain a pair of disjoint 3-GDDs of type
$(24m+14)^t u^1$.

Case 4: $g\equiv 20$ (mod 24), $t=6, 10$,  and $6(t-1)\leq u\leq
7(t-1)$.  Let $g=24k+20$, where $k\geq 1$. Employ a 3-GDD of type
$8^{3k+1} 12^1$ and disjoint pairs of 3-GDDs of types
 $8^t u^1$ and  $12^t u^1$ (Case 1 and Corollary \ref{g6x}). Then obtain a pair of disjoint 3-GDDs of type
$(24k+20)^t u^1$.\qed

\section{The case  $g\equiv 5$ $({\rm mod}$ $6)$}

We shall solve the existence problem of a pair of disjoint modified
group divisible designs in this section. By doing so, the case
$g\equiv 5$ $({\rm mod}$ $6)$ will be completed.

 Let $X$ be a finite set of $gt$
points and $K$ a set of positive integers. A {\em modified group
divisible design} (introduced by Assaf in \cite{assaf}) $K$-GDD is a
quadruple $(X, {\cal G},{\cal H}, {\cal A})$ satisfying the
following properties: $(1)$ $\cal G$ is a partition of $X$ into $t$
$g$-subsets $G_i=\{x_{i,0},x_{i,1},\ldots,x_{i,g-1}\}$, $0\leq i\leq
t-1$. Each $G_i$ is called a {\em group}. ${\cal H}$ is  a partition
of $X$ into $g$ $t$-subsets
$H_j=\{x_{0,j},x_{1,j},\ldots,x_{t-1,j}\}$, $0\leq j\leq g-1$. Each
$H_j$ is called a {\em hole}; $(2)$ $\cal A$ is a set of subsets of
$X$ (called {\em blocks}), each of cardinality from $K$, such that a
block contains no more than one point of any group and any hole;
$(3)$ every pair of points from distinct groups and distinct holes
occurs in exactly one block. A modified group divisible design
$\{3\}$-GDD with $t$ groups and $g$ holes is denoted by
3-MGDD$(g,t)$. Notice that a 3-MGDD$(g,t)$ can also be regarded as a
3-MGDD$(t,g)$. The necessary conditions of the existence of a
3-MGDD$(g,t)$ are $g,t\geq 3$, $(g-1)(t-1)\equiv 0$ (mod 2), and
$gt(g-1)(t-1)\equiv 0$ (mod 6). Similarly, a pair of disjoint
3-MGDD$(g,t)$s means two 3-MGDD$(g,t)$s having same group set and
hole set but disjoint block sets. A 3-MGDD$(3,t)$ is actually same
as an ITD$(3,t)$. So there does not exist a pair of disjoint
3-MGDD(3,3)s. We shall show that it is the only exception.

\begin{Lemma}\label{mgddc}
Suppose that there exists a $(v,K,1)$-PBD. If there exists a pair of
disjoint $3$-MGDD$(g,k)$s for any $k\in K$, then so does a pair of
disjoint $3$-MGDD$(g,v)$s.
\end{Lemma}
\proof Let $(X,{\cal B})$ be a $(v,K,1)$-PBD, ${\cal
G}=\{\{x\}\times I_g:x\in X\}$,  and ${\cal H}=\{X\times\{i\} :i\in
I_g\}$. For any block $B\in{\cal B}$, construct a  pair of disjoint
$3$-MGDD$(g,{|B|})$s  with group set ${\cal G}_B=\{\{x\}\times
I_g:x\in B\}$,  hole set ${\cal H}_B=\{B\times\{i\} :i\in I_g\}$,
and disjoint block sets ${\cal A}_B^1$ and ${\cal A}_B^2$. Define
${\cal A}^1=\cup_{B\in {\cal B}}{\cal A}_B^1$ and ${\cal
A}^2=\cup_{B\in {\cal B}}{\cal A}_B^2$. Then it is immediate that
$(X,{\cal G},{\cal H},{\cal A}^1)$ and $(X,{\cal G},{\cal H},{\cal
A}^2)$ are two disjoint $3$-MGDD$(g,v)$s. \qed

\begin{Lemma}\label{pbd}{\rm (\cite{crc2})}
{\rm(1)} There exists a $(v,\{3,4,6\},1)$-PBD for any  $v\equiv 0,1$
{\rm(mod $3$)}. {\rm(2)}  There exists a $(v,\{3,5\},1)$-PBD for any
$v\equiv 1$ {\rm(mod $2$)}.
\end{Lemma}

{\Lemma \label{mgddsm} For $t=4,6$, there exists a pair of disjoint
$3$-MGDD$(5,t)$s.  }

 \proof
 (1) Let ${\cal G}=\{\{i,i+1,i+2,i+3,i+4\}:i=0,5,10,15\}$ and ${\cal
 H}=\{\{j,j+5,j+10,j+15\}:j=0,1,2,3,4\}$. We construct directly a
 pair of disjoint 3-MGDD$(5,4)$s $(I_{20},{\cal G},{\cal H},{\cal
 A}_1)$ and $(I_{20},{\cal G},{\cal H},{\cal
 A}_2)$, where the blocks are listed below.
\begin{center}
\begin{tabular}{lllllllllll}
${\cal A}_1:$ &$\{0,6,12\}$ &$\{0,7,11\}$ & $\{0,8,16\}$ &
$\{0,9,17\}$ & $\{0,13,19\}$ & $\{0,14,18\}$ \\

& $\{1,5,12\}$ &$\{1,7,18\}$ & $\{1,8,14\}$ & $\{1,9,15\}$
&$\{1,10,19\}$ &
$\{1,13,17\}$\\

& $\{2,5,13\}$ & $\{2,6,19\}$ & $\{2,8,15\}$ & $\{2,9,10\}$ &
$\{2,11,18\}$ & $\{2,14,16\}$ \\

& $\{3,5,16\}$ & $\{3,6,14\}$ &$\{3,7,19\}$ & $\{3,9,11\}$ &
$\{3,10,17\}$ & $\{3,12,15\}$ \\

& $\{4,5,18\}$ & $\{4,6,17\}$ & $\{4,7,13\}$ & $\{4,8,10\}$ &
$\{4,11,15\}$ & $\{4,12,16\}$ \\

&$\{5,11,19\}$ & $\{5,14,17\}$ & $\{6,10,18\}$ & $\{6,13,15\}$ &
$\{7,10,16\}$ & $\{7,14,15\}$ \\

& $\{8,11,17\}$ & $\{8,12,19\}$ & $\{9,12,18\}$ & $\{9,13,16\}$
\medskip\\

${\cal A}_2:$ &$\{0,6,13\}$ & $\{0,7,14\}$ & $\{0,8,17\}$ &
$\{0,9,16\}$ & $\{0,11,18\}$ & $\{0,12,19\}$\\

 & $\{1,5,19\}$ &$\{1,7,15\}$ & $\{1,8,12\}$ & $\{1,9,13\}$ &$\{1,10,17\}$ &
$\{1,14,18\}$ \\

& $\{2,5,18\}$ & $\{2,6,15\}$ & $\{2,8,14\}$ & $\{2,9,11\}$ &
$\{2,10,16\}$ & $\{2,13,19\}$ \\

& $\{3,5,11\}$ & $\{3,6,19\}$ &$\{3,7,10\}$ & $\{3,9,17\}$ &
$\{3,12,16\}$ & $\{3,14,15\}$ \\

& $\{4,5,12\}$ & $\{4,6,10\}$ & $\{4,7,18\}$ & $\{4,8,16\}$ &
$\{4,11,17\}$ & $\{4,13,15\}$\\
 &$\{5,13,17\}$ & $\{5,14,16\}$ &$\{6,12,18\}$ & $\{6,14,17\}$ & $\{7,11,19\}$ & $\{7,13,16\}$ \\
 &
$\{8,10,19\}$ & $\{8,11,15\}$ & $\{9,10,18\}$ & $\{9,12,15\}$
\end{tabular}
\end{center}

 (2) Let $X=(Z_5\times
 I_5)\cup\{\infty_i:i\in I_5\}$, ${\cal G}=\{\{x\}\times I_5:x\in Z_5\}\cup\{\infty_i:i\in I_5\}$, and
 ${\cal H}=\{(Z_5\times\{i\})\cup\{\infty_i\}:i\in I_5\}$. A 3-MGDD$(5,6)$ is constructed
  on $X$ in \cite{assaf} with group set ${\cal G}$, hole set ${\cal H}$ and block sets ${\cal B}_1$
  developed under (mod 5, $-$) by the following blocks:
  \begin{center}
\begin{tabular}{lll}
$\{(0,0),(1,1),(3,2)\}$ & $\{(0,0),(1,2),(2,4)\}$
&$\{(0,1),(3,2),(2,3)\}$\\
 $\{(0,0),(3,1),(1,3)\}$ &
$\{(0,2),(1,3),(4,4)\}$ & $\{(0,1),(1,2),(3,4)\}$\\
$\{(0,0),(2,3),(1,4)\}$ &$\{(0,0),(2,2),(4,3)\}$
&$\{(0,0),(2,1),(3,4)\}$\\
$\{(0,1),(1,3),(2,4)\}$ &$\{(0,0),(4,1),\infty_0\}$&$\{(0,2),(3,3),\infty_0\}$ \\

$\{(0,0),(4,2),\infty_1\}$ & $\{(0,1),(4,4),\infty_1\}$&$\{(0,0),(3,3),\infty_2\}$ \\

 $\{(0,2),(3,4),\infty_2\}$ &
$\{(0,0),(4,4),\infty_3\}$&$\{(0,1),(4,3),\infty_3\}$\\

$\{(0,1),(4,2),\infty_4\}$ & $\{(0,3),(2,4),\infty_4\}$
\end{tabular}
\end{center}
Let ${\cal
B}_2=\{\{(x,a+1),(y,b+1),(z,c+1)\}:\{(x,a),(y,b),(z,c)\}\in{\cal
B}_1\}$, where $\infty_i+1=\infty_i$ for $i\in I_5$. It is readily
checked that ${\cal B}_1$ and ${\cal B}_2$ form block sets of two
disjoint 3-MGDD$(5,6)$s.    \qed

{\Lemma \label{mgdd1} There exists a pair of disjoint
$3$-MGDD$(g,t)$s for any one of the following parameters{\rm :}

\begin{enumerate}
\item[{\rm (1)}] $g\geq 4$ and $t=3;$

\item[{\rm (2)}] $g\equiv 1,3$ {\rm (mod $6$)}, $g\geq 4$ and $t=4,5,6;$

\item[{\rm (3)}] $g\equiv 0,4$ {\rm (mod $6$)}, $g\geq 4$ and $t=5;$

\item[{\rm (4)}] $g\equiv 5$ {\rm (mod $6$)}, $g\geq 5$ and $t=4,6$.
\end{enumerate}
}

 \proof A pair of disjoint 3-MGDD($g$,3)s with $g\geq 4$ exists by
Lemma \ref{itd}.

 For $g\equiv 1,3$ {\rm (mod $6$)}, $g\geq 4$ and $t=4,5,6$, since there are an S$(2,3,g)$ and
  a pair of disjoint 3-MGDD($t$,3)s, we obtain a pair of disjoint 3-MGDD($g$,$t$)s by Lemma \ref{mgddc}.

For $g\equiv 0,4$ {\rm (mod $6$)}, there is a $(g,\{3,4,6\},1)$-PBD
by Lemma \ref{pbd}. A pair of disjoint 3-MGDD($5$,3)s exists by the
above discussion. And a pair of disjoint 3-MGDD($5$,4)s and a pair
of disjoint 3-MGDD($5$,6)s are given in Lemma \ref{mgddsm}. So  we
obtain a pair of disjoint 3-MGDD($g$,5)s by Lemma \ref{mgddc}.

For $g\equiv 5$ {\rm (mod $6$)} and $t=4,6$, there is a
$(g,\{3,5\},1)$-PBD by Lemma \ref{pbd}. Utilize
  pairs of disjoint 3-MGDD($t$,3)s  and disjoint 3-MGDD($t$,5)s. And then obtain a pair of disjoint 3-MGDD$(g,t)$s again by Lemma
\ref{mgddc}. \qed

{\Lemma \label{mgdd} Let $g$ and $t$ be positive integers satisfying
$g,t\geq 3$, $(g,t)\neq (3,3)$, $(g-1)(t-1)\equiv 0$ {\rm (mod $2$)}
and $gt(g-1)(t-1)\equiv 0$ {\rm(mod $6$)}. Then there exists a pair
of disjoint $3$-MGDD$(g,t)$s. }

\proof The conclusion follows by using Lemmas \ref{mgddc}, \ref{pbd}
and \ref{mgdd1}. So we only point out the main ingredients. For
$t\equiv 1,3$ (mod 6), $t\geq 3$ and $g\geq 4$, use an S$(2,3,t)$
and a pair disjoint 3-MGDD$(g,3)$s.  If $t\equiv 2$ (mod 6), then
$t\geq 8$, $g\geq 3$ and $g\equiv 1,3$ (mod 6). Use an S$(2,3,g)$
and a pair disjoint 3-MGDD$(t,3)$s. If $t\equiv 5$ (mod 6), then
$g\equiv 0,1$ (mod 3) and $g\geq 4$. Use a $(t,\{3,5\},1)$-PBD and a
pair of disjoint 3-MGDD$(g,s)$s for $s=3,5$. If $t\equiv 0,4$ (mod
6), then $g\geq 3$ is odd. Use a $(t,\{3,4,6\},1)$-PBD and a pair of
disjoint 3-MGDD$(g,s)$s for $s=3,4,6$. \qed

The following lemmas deal with the admissible triples $(g,t,u)$ with
$g\equiv 5$ (mod 6), so either $u\equiv 1$ (mod 2) if $t\equiv 0$
(mod 6), or  $u\equiv 5$ (mod 6) if $t\equiv 2$ (mod 6), or $u\equiv
3$ (mod 6) if $t\equiv 4$ (mod 6).

 {\Lemma \label{goddt04} Let $(g,t,u)$ be any
admissible triple with $g\equiv 1,5$ {\rm(mod $6$)}, $t\equiv 0,4$
{\rm(mod $6$)}, $g\geq 5$, $t\geq 4$, and $u\leq t-1$. Then there
exists a pair of disjoint $3$-GDDs of type $g^tu^1$. }

\proof For  $g\equiv 1,5$ {\rm(mod $6$)}, $t\equiv 0,4$ {\rm(mod
$6$)},  $g\geq 5$, and $t\geq 4$, by Lemma \ref{mgdd} there is a
pair of disjoint 3-MGDD$(g,t)$s on a $gt$-set $X$ with  group set
${\cal G}$, hole set ${\cal H}$ and disjoint block sets ${\cal A}_1$
and ${\cal A}_2$. Further $(1,t,u)$ is also an admissible triple.
Let $U$ be a $u$-set disjoint with $X$. For each $H\in {\cal H}$,
construct on $H\cup U$ a pair of disjoint 3-GDDs of type $1^t u^1$
with $U$ as the long group and ${\cal B}_H^1$ and ${\cal B}_H^2$ as
the block sets. For $i=1,2$, let ${\cal C}_i={\cal
A}_i\cup(\cup_{H\in{\cal H}}{\cal B}_H^i)$. Thus $(X,{\cal
G}\cup\{U\},{\cal C}_1)$ and $(X,{\cal G}\cup\{U\},{\cal C}_2)$ form
a pair of disjoint $3$-GDDs of type $g^tu^1$. \qed

{\Lemma \label{511sm}  There exists a pair of disjoint $3$-GDDs of
type $g^tu^1$, where
$(g,t,u)\in\{(5,4,3),(5,4,9),(11,4,3),(11,4,9),(11,4,15),(11,4,21),(11,4,27),(11,8,5),(11,6,7),$
 \\ $(11,6,9)\}$.
}

\proof For
$(g,t,u)=(5,4,3),(5,4,9),(11,4,3),(11,4,9),(11,4,15),(11,4,21),(11,4,$
$27),(11,8,5)$, let $D=\{1,2,\ldots,$
$gt/2\}\setminus\{t,2t,\ldots,[g/2]t\}$. Since a partition of $D$
into $D_1$ and $D_2$ satisfying the conditions of Lemma \ref{dif} is
given in Section 5 of \cite{chr}, there exists a pair of disjoint
$3$-GDDs of type $g^tu^1$. For $g=11,t=6$ and $u=7,9$, apply the
Weighting Construction to a TD$(7,7)$ as in \cite[Lemma 5.4]{chr}.
Take a block of the TD(7,7) and weight 5 to six points and weight 1
or 3 to the other point of the block. Then weight 1 to all the other
points. Since there is a pair of disjoint 3-GDDs of type $1^7$,
$1^63^1$, $1^65^1$, or $5^63^1$ (Lemmas \ref{g1} and \ref{goddt04}),
a pair of disjoint $3$-GDDs of type $g^tu^1$ also exists. \qed

{\Lemma \label{goddt2} Let $(g,t,u)$ be any admissible triple with
$g=5,11$, $u< g$, $t\equiv 2$ {\rm(mod $6$)}, and $t\geq 14$. Then
there exists a pair of disjoint $3$-GDDs of type $g^tu^1$.}

\proof For $g=5,11$, $u< g$, $t\equiv 2$ {\rm(mod $6$)}, and $t\geq
14$, there is a pair of disjoint 3-GDDs of type
$(2g)^{(t-6)/2}(5g+u)^1$ by Lemma \ref{g6k+4}. There exists a pair
of disjoint 3-GDDs of type $g^{t-6}(6g+u)^1$ by Corollary
\ref{2gtog}. There exists  a pair of disjoint 3-GDDs of type
$g^{6}u^1$ by Lemmas \ref{goddt04}
 and \ref{511sm}. So a pair of disjoint $3$-GDDs of type
$g^tu^1$ exists by Filling Construction II. \qed

{\Lemma \label{g511} The Main Theorem holds for any admissible
triple $(g,t,u)$  with $g\equiv 5$ {\rm (mod $6$)} and $5\leq
g\leq29$.}

\proof The case of $u>4(t-1)$ is solved by Lemma \ref{most}.  Also
noting that for  $t\geq 6$ (must be even) and $u>g$, there exists a
pair of disjoint $3$-GDDs of type $g^tu^1$ by Corollary \ref{2gtog}
since there is a pair of disjoint $3$-GDDs of type
$(2g)^{t/2}(u-g)^1$ by Lemma \ref{g6k+4}, we only need to consider
the cases $u\leq 12$ if $t=4$ and $u\leq4(t-1)$ and $u<g$ if $t\geq
6$. All the possibilities are exhausted as follows (with $(g,t,u)$
admissible):

Case 1:  $g=5,11$, and  $u\leq 4(t-1)$, further $u<g$ if $t\geq 6$.
There are several subcases of $t$. (i) $t=4$. There is a pair of
disjoint 3-GDDs of type $g^tu^1$ by Lemma \ref{511sm}. (ii)
  $t\equiv 2$ (mod 6). If $t=8$, then we use Lemma \ref{511sm} to
  deal with the only possible triple $(11,8,5)$. Otherwise $t\geq
  14$ and Lemma \ref{goddt2} gives the solution. (iii) $t\equiv 0,4$ (mod
  6). If $u\leq t-1$, then we use Lemma \ref{goddt04} to obtain the
  desired pair of 3-GDDs. Otherwise $t-1<u<g$.  Thus all the possible
  admissible  triples are $(11,6,7)$ and $(11,6,9)$, the solutions of which are
  listed in Lemma \ref{511sm}.

Case 2: $g=17$, and $u\leq 4(t-1)$, and further $u<g$ if $t\geq 6$.
Since $(g,t,u)$ is admissible, it is readily checked that $u\leq
3(t-1)$. Hence a pair of disjoint 3-GDDs of type $g^tu^1$ exists by
Construction \ref{glarge} since a 3-GDD of type $3^45^1$ and
disjoint pairs of 3-GDDs of types $3^tu^1$ and $5^tu^1$ exist.

Case 3: $g=29$, and $u\leq 4(t-1)$. Then a pair of disjoint 3-GDDs
of type $g^tu^1$ exists by Construction \ref{glarge} since a 3-GDD
of type $5^49^1$ and disjoint pairs of 3-GDDs of types $5^tu^1$ and
$9^tu^1$ exist.

Case 4: $g=23$, $u\leq4(t-1)$, and $u<g$. If $u\leq 3(t-1)$,  a pair
of disjoint 3-GDDs of type $g^tu^1$ exists by Construction
\ref{glarge} since  a 3-GDD of type $3^65^1$ and disjoint pairs of
 3-GDDs of types $3^tu^1$ and $5^tu^1$ exist. Thus it remains
only to deal with the cases $t=6$ and odd $u$ with $15< u\leq 20$.

Similar to \cite[Lemma 4.3]{chr}, start from a $\{2,3\}$-GDD of type
$1^{18} 5^1$ $(X,{\cal G},{\cal B})$, where $G\in{\cal G}$, $|G|=5$,
and the blocks of size 2 form four parallel classes of $X\setminus
G$, say ${\cal P}_i,i\in I_4$. Let
$U=\{\infty_1,\infty_2,\ldots,\infty_{u}\}$, $Y=(X\times I_6)\cup
U$, and ${\cal H}=\{X\times\{i\}:i\in I_6\}\cup\{U\}$. First for
each $B\in{\cal B}$ and $|B|=3$, construct on $B\times I_6$ a pair
of disjoint ITD$(3,6)$s omitting the idempotent parallel class,
whose group set is $\{\{x\}\times I_6:x\in B\}$ and two block sets
are ${\cal A}_B^1$ and ${\cal A}_B^2$. Then we deal with $G$, the
group of size $5$  in ${\cal G}$. Construct on $(G\times I_6) \cup
U$ a pair of disjoint 3-GDDs of type $5^6u^1$ with group set
$\{G\times \{i\}:i\in I_6\}\cup\{U\}$ and block sets ${\cal D}^1$
and ${\cal D}^2$.
 After that let
$U_k=\{\infty_{5k+1},\infty_{5k+2},\ldots,\infty_{5k+5}\}$, where
$k=0,1,2$, and $U_3=U\setminus(U_0\cup U_1\cup U_2)$. For each pair
$P\in {\cal P}_3$ construct on $(P\times I_6)\cup U_3$ a pair of
disjoint 3-GDDs of type $2^6(u-15)^1$, whose group set is
$\{\{x\}\times I_6:x\in B\}\cup \{U_3\}$ and two block sets are
${\cal E}_P^1$ and ${\cal E}_P^2$. Finally for each ${\cal P}_k$,
$k=0,1,2$,  the set $\{\{(x,i),(y,j)\}:\{x,y\}\in {\cal P}_k, i\neq
j\in I_6\}$  can
 be partitioned  into 5 disjoint 1-factors of $X\setminus
I_6$, denoted by $F_{k0},F_{k1},\ldots,F_{k4}$. Let ${\cal
F}_k^1=\cup_{0\leq l\leq
4}\{\{\infty_{2k+1+l},\alpha,\beta\}:\{\alpha,\beta\}\in F_{kl}\}$
and ${\cal F}_k^2=\cup_{0\leq l\leq
4}\{\{\infty_{2k+1+l},\alpha,\beta\}:\{\alpha,\beta\}\in
F_{k,l+1}\}$.  For $i=1,2$, let ${\cal C}^i={\cal
D}^i\cup(\cup_{B\in{\cal B},|B|=3}{\cal A}_B^i)\cup(\cup_{P\in{\cal
P}_3}{\cal E}_P^i)\cup(\cup_{0\leq k\leq 2}{\cal F}_{k}^i)$. It can
be checked that  $(Y,{\cal H},{\cal C}^1)$ and $(Y,{\cal H},{\cal
C}^2)$ form two disjoint 3-GDDs of type $23^6u^1$. \qed

{\Lemma \label{g6k+5} The Main Theorem holds for any admissible
triple  $(g,t,u)$  with $g\equiv 5$ {\rm (mod $6$)}.}

\proof We can employ Lemma \ref{most} to treat $u>4(t-1)$, Corollary
\ref{t3} to treat $t=3$, and Lemma \ref{g511} to treat $g\leq 29$.
So let $g=6n+5$, $n\geq 5$,  $t\geq 4$ and $u\leq 4(t-1)$. Apply
induction on $n$. Suppose that there is a pair of 3-GDDs of type
$h^s v^1$ for any admissible triple $(h,s,v)$ with $h=6l+5$, and
$l<n$. If $n\equiv 3,5$ (mod 6), then a  3-GDD of type $n^65^1$
exists by Lemma \ref{2.3}. And
 disjoint pairs of 3-GDDs of types $n^tu^1$ and $5^t u^1$ also exist by Lemma \ref{g511} or by
the assumption. So a pair of disjoint 3-GDDs of type $(6n+5)^tu^1$
exists by Construction \ref{glarge}. If $n\equiv 0,4$ (mod 6), or
$n\equiv 1$ (mod 6), or $n\equiv 2$ (mod 6), also utilize
Construction \ref{glarge} but taking instead a 3-GDD of type
$(n-1)^611^1$, or $(n-2)^617^1$, or $(n-3)^623^1$, and so on. This
completes the proof. \qed

\section{Conclusion}
 Summing up the results of  Lemmas  \ref{intersect}, \ref{g6k+3},
\ref{most}, \ref{g6k+4}, \ref{g6k+2}, \ref{g6k+5}, and Corollary
\ref{g6x}, we obtain the Main Theorem.

To end this paper we mention a byproduct on group divisible codes,
which  play an important role in the determination of some optimal
constant-weight and constant-composition codes. Here we do not dwell
on relevant notations on coding theory and the interested readers
are referred to \cite{chee2,zhangge}. If $(X,{\cal G},{\cal B}_1)$
and $(X,{\cal G},{\cal B}_2)$ are a pair of disjoint 3-GDDs of type
$g^tu^1$, from which we can naturally obtain a pair of disjoint
$(n,4,3)_2$ codes ${\cal C}_1$ and ${\cal C}_2$ where $n=gt+u$. As
in \cite{chee1}, replace each occurrence of 1 with $i$ in each
codeword of ${\cal C}_i$ to yield a new code ${\cal C}_i'$
($i=1,2$). Thus ${\cal C}_1'\cup{\cal C}_2'$ forms a ternary group
divisible codes of weight three, distance four and  size $2b$, where
$b=\frac{1}{6}(g^2t(t-1)+2gtu)$, the number of blocks in a 3-GDD of
type $g^tu^1$.

\vspace{1.0cm}\par \centerline{\Large\bf Acknowlegements}
\vspace{0.2cm} A portion of this research was carried out while the
first author was visiting Nanyang Technological University in 2008,
and he wishes to express many thanks to the Division of Mathematics
for their services.

{\bf Appendix}

We list a pair of disjoint 3-GDDs of type $g^t u^1$, where
$(g,t,u)\in\{(3,4,1),(3,4,5),$ $(3,4,7)\}$. The point set is
$I_{gt+u}$. The groups are $\{it+j:i\in I_g\}$, $j\in I_t$, and
$\{gt,gt+1,\ldots,gt+u-1\}$. And the disjoint block sets ${\cal
A}_1$ and ${\cal A}_2$ are as follows.

(1) $(g,t,u)=(3,4,1)$.

\begin{center}
\begin{tabular}{llllllll}
${\cal A}_1:$&$\{12,0,1\}$ & $\{12,2,3\}$ & $\{12,4,6\}$ &
$\{12,5,7\}$ & $\{12,8,11\}$ & $\{12,9,10\}$\\
&$ \{0,2,5\}$ & $\{0,3,6\}$ & $\{0,7,9\}$ & $\{0,10,11\}$ &
$\{1,2,8\}$ & $\{1,3,10\}$\\
&$ \{1,4,11\}$ & $\{1,6,7\}$ & $\{2,4,7\}$ & $\{2,9,11\}$ &
$\{3,4,9\}$ & $\{3,5,8\}$\\
&$ \{4,5,10\}$ & $\{5,6,11\}$ & $\{6,8,9\}$ & $\{7,8,10\}$

\medskip\\

${\cal A}_2:$& $\{12,0,2\}$ & $\{12,1,3\}$ & $\{12,4,5\}$ &
$\{12,6,9\}$ & $\{12,7,8\}$ & $\{12,10,11\}$\\
&$ \{0,1,6\}$ & $\{0,3,5\}$ & $\{0,7,10\}$ & $\{0,9,11\}$ &
$\{1,2,7\}$ & $\{1,4,10\}$\\
&$ \{1,8,11\}$ & $\{2,3,4\}$ & $\{2,5,11\}$ & $\{2,8,9\}$ &
$\{3,6,8\}$ & $\{3,9,10\}$\\
&$ \{4,6,11\}$ & $\{4,7,9\}$ & $\{5,6,7\}$ & $\{5,8,10\}$
\end{tabular}
\end{center}

(2) $(g,t,u)=(3,4,5)$.
\begin{center}
\begin{tabular}{llllllll}
${\cal A}_1:$& $\{12,0,1\} $& $\{12,2,3\} $& $\{12,4,5\} $&
$\{12,6,7\} $& $\{12,8,9\} $& $\{12,10,11\} $ \\

& $\{13,0,2\} $& $\{13,1,3\} $& $\{13,4,6\} $& $\{13,5,7\} $&
$\{13,8,10\} $& $\{13,9,11\} $\\

 & $\{14,0,3\} $& $\{14,1,2\} $&$\{14,4,7\} $& $\{14,5,6\} $& $\{14,8,11\} $& $\{14,9,10\} $ \\

 &$\{15,0,5\} $& $\{15,1,6\} $& $\{15,2,8\} $& $\{15,3,9\} $&
$\{15,4,11\} $& $\{15,7,10\} $ \\
& $\{16,0,10\} $& $\{16,1,11\} $&$\{16,2,7\} $& $\{16,3,4\} $&
$\{16,5,8\} $& $\{16,6,9\} $ \\
& $\{0,6,11\} $& $\{0,7,9\} $& $\{1,4,10\} $& $\{1,7,8\} $&
$\{2,4,9\} $& $\{2,5,11\} $\\
& $\{3,5,10\} $& $\{3,6,8\} $

\medskip\\

${\cal A}_2:$& $\{12,0,2\} $& $\{12,1,3\} $& $\{12,4,6\} $& $\{12,5,7\} $& $\{12,8,10\} $& $\{12,9,11\} $\\
& $\{13,0,1\} $& $\{13,2,3\} $& $\{13,4,5\} $& $\{13,6,7\} $& $\{13,8,9\} $& $\{13,10,11\} $\\
& $\{14,0,5\} $& $\{14,1,4\} $& $\{14,2,8\} $& $\{14,3,10\} $& $\{14,6,11\} $& $\{14,7,9\} $\\
& $\{15,0,11\} $& $\{15,1,10\} $& $\{15,2,5\} $& $\{15,3,4\} $& $\{15,6,9\} $& $\{15,7,8\} $\\
& $\{16,0,6\} $& $\{16,1,7\} $& $\{16,2,9\} $& $\{16,3,8\} $& $\{16,4,11\} $& $\{16,5,10\} $\\
& $\{0,3,9\} $& $\{0,7,10\} $& $\{1,2,11\} $& $\{1,6,8\} $& $\{2,4,7\} $& $\{3,5,6\} $\\
& $\{4,9,10\} $& $\{5,8,11\} $
\end{tabular}
\end{center}

(3) $(g,t,u)=(3,4,7)$.
\begin{center}
\begin{tabular}{llllllll}
${\cal A}_1:$& $\{12,0,1\} $& $\{12,2,3\} $& $\{12,4,5\} $& $\{12,6,7\} $& $\{12,8,9\} $& $\{12,10,11\} $\\
& $\{13,0,2\} $& $\{13,1,3\} $& $\{13,4,6\} $& $\{13,5,7\} $& $\{13,8,10\} $& $\{13,9,11\} $\\
& $\{14,0,3\} $& $\{14,1,2\} $& $\{14,4,7\} $& $\{14,5,6\} $& $\{14,8,11\} $& $\{14,9,10\} $\\
& $\{15,0,5\} $& $\{15,1,4\} $& $\{15,2,8\} $& $\{15,3,9\} $& $\{15,6,11\} $& $\{15,7,10\} $\\
& $\{16,0,6\} $& $\{16,1,7\} $& $\{16,2,9\} $& $\{16,3,8\} $& $\{16,4,10\} $& $\{16,5,11\} $\\
& $\{17,0,10\} $& $\{17,1,11\} $& $\{17,2,5\} $& $\{17,3,4\} $& $\{17,6,9\} $& $\{17,7,8\} $\\
& $\{18,0,11\} $& $\{18,1,10\} $& $\{18,2,7\} $& $\{18,3,6\} $& $\{18,4,9\} $& $\{18,5,8\} $\\
& $\{0,7,9\} $& $\{1,6,8\} $& $\{2,4,11\} $& $\{3,5,10\} $
\medskip\\
${\cal A}_2:$ & $\{12,0,2\} $& $\{12,1,3\} $& $\{12,4,6\} $& $\{12,5,7\} $& $\{12,8,10\} $& $\{12,9,11\} $\\
& $\{13,0,1\} $& $\{13,2,3\} $& $\{13,4,5\} $& $\{13,6,7\} $& $\{13,8,9\} $& $\{13,10,11\} $\\
& $\{14,0,5\} $& $\{14,1,4\} $& $\{14,2,8\} $& $\{14,3,9\} $& $\{14,6,11\} $& $\{14,7,10\} $\\
& $\{15,0,3\} $& $\{15,1,2\} $& $\{15,4,7\} $& $\{15,5,6\} $& $\{15,8,11\} $& $\{15,9,10\} $\\
& $\{16,0,10\} $& $\{16,1,11\} $& $\{16,2,4\} $& $\{16,3,5\} $& $\{16,6,8\} $& $\{16,7,9\} $\\
& $\{17,0,11\} $& $\{17,1,10\} $& $\{17,2,7\} $& $\{17,3,6\} $& $\{17,4,9\} $& $\{17,5,8\} $\\
& $\{18,0,7\} $& $\{18,1,6\} $& $\{18,2,9\} $& $\{18,3,8\} $& $\{18,4,11\} $& $\{18,5,10\} $\\
& $\{0,6,9\} $& $\{1,7,8\} $& $\{2,5,11\} $& $\{3,4,10\} $
\end{tabular}
\end{center}

\end{document}